\documentclass[11pt,a4paper,reqno]{amsart}
\usepackage[T1]{fontenc}
\usepackage{lmodern}
\usepackage[margin=26mm]{geometry}
\usepackage{amsmath,amssymb,amsthm,mathtools}
\usepackage{booktabs}
\usepackage{microtype}
\usepackage[hidelinks]{hyperref}
\hypersetup{
 pdftitle={Polynomial approximations to Calabi-Yau metrics on quartic K3 surfaces and quintic threefolds},
 pdfauthor={Kefeng Liu and Liangdi Zhang}
}
\usepackage{enumitem}
\setlist[enumerate]{itemsep=3pt,topsep=5pt}
\numberwithin{equation}{section}

\newtheorem{theorem}{Theorem}[section]

\newtheorem{maintheorem}{Theorem}

\newtheorem{lemma}[theorem]{Lemma}
\newtheorem{proposition}[theorem]{Proposition}
\newtheorem{corollary}[theorem]{Corollary}
\theoremstyle{remark}
\newtheorem{remark}[theorem]{Remark}
\theoremstyle{plain}

\newcommand{\Q}{\mathbb Q}
\newcommand{\PP}{\mathbb P}
\newcommand{\ddbar}{\sqrt{-1}\partial\bar\partial}

\newcommand{\HS}{\mathrm{HS}}
\newcommand{\op}{\mathrm{op}}

\title[Polynomial Calabi--Yau approximations]{Polynomial approximations to Calabi--Yau metrics on quartic K3 surfaces and quintic threefolds}
\author{Kefeng Liu}
\address{Mathematical Sciences Research Center, Chongqing University of
Technology, Chongqing 400054, China}
\email{kefengliu@cqut.edu.cn}
\author{Liangdi Zhang}
\address{Mathematical Sciences Research Center, Chongqing University of
Technology, Chongqing 400054, China}
\email{ldzhang91@163.com}
\date{}
\makeatletter
\@ifundefined{subjclassname@2020}{%
  \@namedef{subjclassname@2020}{\textup{2020} Mathematics Subject Classification}%
}{}
\makeatother
\begin{document}
\begin{abstract}
On the Fermat quartic K3 surface and the Fermat quintic
threefold, we construct K\"ahler metrics from explicit rational
polynomial potentials and bound their differences from the
Ricci-flat metrics in integral tensor norms. For specified compact
families of smooth quartic and quintic hypersurfaces near the
Fermat hypersurfaces, we give a terminating algorithm that
constructs polynomial K\"ahler metrics with prescribed pointwise
bounds on their differences from the Ricci-flat metrics and on
their Ricci forms. The resulting metrics converge to the
Ricci-flat metrics in every $C^r$ norm, uniformly in the parameter.
\end{abstract}
\maketitle

\section{Introduction}
In Problem~23 of his 1982 problem list, Yau
\cite[p.~676]{Yau82} asks: \textquotedblleft Construct an explicit
metric with zero Ricci curvature on the K3 surface.\textquotedblright{}
The existence of such metrics follows from Yau's solution
\cite{Yau} of Calabi's conjecture \cite{Calabi}: on a compact
K\"ahler manifold with vanishing real first Chern class, each
K\"ahler class contains a unique Ricci-flat K\"ahler metric.
The theorem does not give an explicit formula. Here we seek
explicit approximations in the hyperplane class on quartic K3
surfaces and quintic Calabi--Yau threefolds, whose defining
equations and holomorphic volume forms are explicit; see, for
example, Huybrechts \cite[Chapter~1]{Huybrechts} and Joyce
\cite{Joyce}.

Following a proposal of Yau, Tian \cite{Tian} proved approximation
of polarized K\"ahler metrics by projectively induced metrics.
For Bergman kernel asymptotics and their coefficients, see,
for example, Catlin \cite{Catlin}, Zelditch \cite{Zelditch},
Lu \cite{Lu}, and Dai, Liu and Ma \cite{DaiLiuMa}.
Ma and Marinescu \cite{MaMarinescu} give a systematic treatment.
Cao and Keller \cite[Theorem~1.1 and Corollary~4.2]{CaoKeller}
prove smooth convergence of normalized volume-balanced metrics
to the solution of the prescribed-volume Monge--Amp\`ere equation.

Gluing gives approximations in degenerating K\"ahler classes.
Gross and Wilson \cite{GrossWilson} construct approximating
metrics on elliptically fibered K3 surfaces in a collapsing
limit. Fu, Liu, Xiao and Yau \cite[Theorem~1.1]{FLXY} give
a convergent recursive representation on a fixed Kummer K3
surface in K\"ahler classes approaching the orbifold limit.

Explicit metrics are also needed in string compactification.
Candelas, Horowitz, Strominger and Witten \cite{CHSW} established
the role of Calabi--Yau threefolds in superstring vacua.
The normalization of physical Yukawa couplings depends on
metric inner products, which are not determined by topology;
see, for example, Mirjani\'c and Mishra \cite[Section~2.2]{MM}.

For numerical constructions on compact manifolds, see, for
example, Headrick and Wiseman \cite{HeadrickWiseman}, who
compute Ricci-flat metrics on Kummer surfaces. Donaldson \cite{Donaldson} develops
numerical projective methods, and Douglas, Karp, Lukic and
Reinbacher \cite{DKLR} apply balanced metrics to a one-parameter
family of quintics. Braun, Brelidze, Douglas and Ovrut
\cite{Braun} treat quotients and complete intersections.
Headrick and Nassar \cite{HN} study energy functionals for
approximate metrics.

For neural-network constructions, see,
for example, Anderson, Gerdes, Gray, Krippendorf, Raghuram
and Ruehle \cite{Anderson}, Douglas, Lakshminarasimhan and Qi
\cite{DLQ}, and Larfors, Lukas, Ruehle and Schneider
\cite{LLRS}.

Some of these computations lead to explicit potentials.
Mirjani\'c and Mishra \cite[Section~5]{MM} use symmetries of
ambient extensions to obtain compact symbolic approximations.
Lee and Lukas \cite{LL} fit learned potentials to Donaldson's
ansatz and obtain explicit parameter-dependent approximations
on quintic and bicubic families.

To measure the accuracy of an explicit potential, one needs
bounds for the difference between its metric and the Ricci-flat
metric. Douglas, Platt, Qi and Barbosa
\cite[Proposition~3.1 and Corollary~3.5]{DPQB} prove such bounds
when the volume-density difference is small in a H\"older norm
involving first derivatives.

We give polynomial potentials with rational coefficients and
prove global bounds for their metrics.
Theorem~\ref{thm:explicit-main} specifies every coefficient of
two potentials on the Fermat quartic and quintic, proves that
they define K\"ahler metrics, and gives tensor $L^1$ bounds
relative to the Ricci-flat metrics.
Theorem~\ref{thm:strong-intro} gives a polynomial construction
with arbitrarily high pointwise accuracy, uniformly on compact
families near these hypersurfaces; the metrics converge in
every $C^r$ norm.
Neither result gives a finite formula for the exact Ricci-flat metric.

For $d=4,5$, put $n=d-2$ and
\[
 X_d=\left\{\sum_{i=0}^{d-1}Z_i^d=0\right\}\subset\PP^{d-1},
 \qquad \omega_0=\omega_{\mathrm{FS}}|_{X_d}.
\]
Here $\omega_{\mathrm{FS}}=\sqrt{-1}\partial\bar\partial\log\sum_i|Z_i|^2$.
Let $\beta_d\in[\omega_0]$ be the Ricci-flat form. For a K\"ahler
form $\gamma$, set $d\mu_\gamma=\gamma^n/\int_{X_d}\gamma^n$.
If $\lambda_i$ are the eigenvalues of a real $(1,1)$ form $\eta$
with respect to $\gamma$, put
$|\eta|_{\gamma,\HS}=(\sum_i\lambda_i^2)^{1/2}$.
Set
\[
 s_j([Z])=\frac{\sum_i|Z_i|^{2j}}
                   {(\sum_i|Z_i|^2)^j}.
\]

For $b=(b_2,\ldots,b_d)$, define
\begin{equation}
 P_A=\frac1{D_A}\sum_b a_b\prod_{j=2}^d s_j^{b_j},\qquad
 \begin{array}{c|cc}
 A&d&D_A\\ \hline
 Q_{12}&4&640000000000\\
 T_{12}&5&7680000000000.
 \end{array}
 \label{eq:recovered-Q12-T12}
\end{equation}
The integers $a_b$ are listed in Tables~\ref{tab:recovered-Q12}
and~\ref{tab:recovered-T12} in Appendix~\ref{app:coefficients};
every unlisted coefficient is zero.

\begin{maintheorem}[Explicit polynomial metrics]
\label{thm:explicit-main}
Let $P_4=P_{Q_{12}}$ and $P_5=P_{T_{12}}$ be defined by
\eqref{eq:recovered-Q12-T12}. Then
\[
 \widehat\omega_d=\omega_0+\sqrt{-1}\partial\bar\partial P_d>0
 \quad(d=4,5),
\]
and
\begin{equation}
 \max_{\gamma\in\{\widehat\omega_d,\beta_d\}}
 \int_{X_d}|\widehat\omega_d-\beta_d|_{\gamma,\HS}
                                      \,d\mu_\gamma
 <\begin{cases}
     3/1000,&d=4,\\
     31/1000,&d=5.
   \end{cases}
 \label{eq:explicit-main-L1}
\end{equation}
\end{maintheorem}

The potentials have $33$ terms on the quartic and $46$ on
the quintic. Lemma~\ref{lem:certified-density} verifies their
global volume-density bounds by finite integer arithmetic;
Lemma~\ref{lem:volume-contrast} gives the tensor estimates.
A shorter quartic potential with $14$ terms is given in
\eqref{eq:recovered-Q8}.

For approximation with arbitrarily small bounds, consider
\begin{equation}
 X_t=\left\{
 F_t(Z):=\sum_{i=0}^{d-1}Z_i^d+
                  \sum_{|\gamma|=d}t_\gamma Z^\gamma=0
 \right\}\subset\PP^{d-1}.
 \label{eq:family}
\end{equation}
Put
\begin{equation}
 |t|_1=\sum_{|\gamma|=d}|t_\gamma|,\qquad
 B_4=\{|t|_1\le1/8\},\qquad B_5=\{|t|_1\le1/16\}.
 \label{eq:parameter-balls}
\end{equation}
Let $\omega_{0,t}=\omega_{\mathrm{FS}}|_{X_t}$, let $g_{0,t}$
be its Riemannian metric, and let $\beta_t$ be the Ricci-flat
form in $[\omega_{0,t}]$. Write
\[
 \Pi_{ab}([Z])=\frac{Z_a\overline Z_b}{\sum_i|Z_i|^2}.
\]
Polynomials in $(t,\Pi)$ are understood as polynomials in
their real and imaginary parts. Derivatives and $C^r$ norms
below are taken along the fibers with respect to $g_{0,t}$.

\begin{maintheorem}[Uniform smooth approximation]
\label{thm:strong-intro}
For $d\in\{4,5\}$, there is an algorithm which, for every
integer $j\ge1$, returns a polynomial $P_j(t,\Pi)$ with rational
coefficients after finitely many steps. For $B_d$ in
\eqref{eq:parameter-balls}, the forms
\[
 \alpha_{j,t}=\omega_{0,t}+\sqrt{-1}\partial_{X_t}\bar\partial_{X_t}
                                   P_j(t,\Pi)\quad(t\in B_d)
\]
are K\"ahler and satisfy
\begin{equation}
 \sup_{t\in B_d}\max\left\{
 \|\alpha_{j,t}-\beta_t\|_{C^0(\omega_{0,t}),\op},
 \|\operatorname{Ric}(\alpha_{j,t})\|_{C^0(\omega_{0,t}),\op}
 \right\}<2^{-j}.
 \label{eq:intro-effective-uniform}
\end{equation}
For every integer $r\ge0$, the sequence $(\alpha_{j,t})_{j\ge1}$
returned by this algorithm satisfies
\begin{equation}
 \sup_{t\in B_d}\|\alpha_{j,t}-\beta_t\|_{C^r(g_{0,t})}
       \le C_{d,r}2^{-j}\qquad(j\ge J_{d,r}),
 \label{eq:intro-smooth-approximation}
\end{equation}
where $C_{d,r}<\infty$ and the integer $J_{d,r}\ge1$ depend
only on $d$ and $r$.
\end{maintheorem}

The algorithm enumerates polynomials with rational coefficients
and tests explicit polynomial inequalities on the family
\eqref{eq:family}. For each $j$, the search stops at the first
polynomial satisfying the required inequalities. The proof gives
no bound on its degree or on the number of steps, and does
not assert optimality.

The constants $C_{d,r}$ and $J_{d,r}$ are not evaluated;
the search does not use them. The algorithm in
Theorem~\ref{thm:strong-intro} has not been run to compute
additional potentials.

Section~\ref{sec:preliminaries} fixes the notation.
Section~\ref{sec:explicit-potential} proves
Theorem~\ref{thm:explicit-main}.
Sections~\ref{sec:projective-background} and~\ref{sec:pointwise}
establish the estimates needed for
Theorem~\ref{thm:strong-intro}, which is proved in
Section~\ref{subsec:parametric-uniform}. That section also
treats a single hypersurface and invariant potentials.
Appendix~\ref{app:coefficients} lists the coefficients.

\section{Preliminaries}\label{sec:preliminaries}
The Fubini--Study form is
$\omega_{\mathrm{FS}}=\sqrt{-1}\partial\bar\partial
\log\sum_i|Z_i|^2$. If
$\omega=\sqrt{-1}g_{i\bar j}\,dz_i\wedge d\bar z_j$,
the associated real metric is
$2\operatorname{Re}(g_{i\bar j}\,dz_i\otimes d\bar z_j)$.
We use
\[
 \Delta_\omega v=g^{i\bar j}\partial_i\partial_{\bar j}v,
 \qquad d\mu_\omega=\frac{\omega^n}{\int_X\omega^n}.
\]
Thus $\Delta_{\mathbb R}=2\Delta_\omega$ and
$|\nabla v|^2=2|\partial v|_\omega^2$.
For the hypersurface $X_t$ in \eqref{eq:family}, set
\[
 d\mu_{0,t}=\frac{\omega_{0,t}^n}{\int_{X_t}\omega_{0,t}^n},
 \qquad d\nu_t=\frac{\beta_t^n}{\int_{X_t}\beta_t^n}.
\]
Both measures have total volume one.

Let $\eta$ be a real $(1,1)$ form, and let
$\lambda_1,\ldots,\lambda_n$ be the eigenvalues of the
Hermitian endomorphism $\omega^{-1}\eta$. We define
\[
 |\eta|_{\omega,\op}=\max_i|\lambda_i|,\qquad
 |\eta|_{\omega,\HS}=\left(\sum_i\lambda_i^2\right)^{1/2},
 \qquad
 \|\eta\|_{C^0(\omega),\op}=\sup_X|\eta|_{\omega,\op}.
\]
For a matrix $A$, $\|A\|_{\op}=\sup_{v\ne0}|Av|/|v|$,
where $|\cdot|$ is the Euclidean norm.
The real symmetric tensor associated to $\eta$ has
Hilbert--Schmidt norm $\sqrt{2}|\eta|_{\omega,\HS}$.

Write $g_{0,t}$ for the real metric of $\omega_{0,t}$ and
$\nabla_{0,t}$ for its Levi-Civita connection. For a real
covariant tensor $T$ on $X_t$, set
\[
 \|T\|_{C^r(g_{0,t})}
 =\max_{0\le k\le r}\sup_{X_t}
       |\nabla_{0,t}^{\,k}T|_{g_{0,t}},
\]
using real Hilbert--Schmidt norms. All these derivatives are
taken on $X_t$, with $t$ fixed.

On a compact connected K\"ahler manifold $(X,\omega)$ of
complex dimension $n\ge2$, choose constants $S\ge1$ and
$\lambda>0$, depending only on $(X,\omega)$, such that the
Poincar\'e inequality for smooth real functions $w$ with
$\int_Xw\,d\mu_\omega=0$,
\begin{equation}
 \|w\|_2^2\le\lambda^{-1}\|\partial w\|_2^2,
 \label{eq:background-P}
\end{equation}
and the Sobolev inequality for every smooth real function $w$,
\begin{equation}
 \|w\|_{2n/(n-1)}^2
 \le S\bigl(\|\partial w\|_2^2+\|w\|_2^2\bigr),
 \label{eq:background-S}
\end{equation}
hold; see, for example, Li \cite[Chapters~9 and~14]{Li}.
Here $\|w\|_p=(\int_X|w|^p\,d\mu_\omega)^{1/p}$ and
$\|\partial w\|_2^2=\int_X|\partial w|_\omega^2\,d\mu_\omega$.

\begin{lemma}\label{lem:family-smoothness}
The hypersurfaces $X_t$ in \eqref{eq:family} are smooth and
connected for all $t$ in an open neighborhood of the set $B_d$
defined in \eqref{eq:parameter-balls}. Their union is a smooth
family over that neighborhood.
\end{lemma}
\begin{proof}
Write $F_t=F_0+G_t$ and normalize $|Z|=1$. The vector of
complex derivatives satisfies
\[
 |\nabla F_0|^2=d^2\sum_i|Z_i|^{2d-2}\ge d^{4-d},
 \qquad
 |\nabla Z^\gamma|\le
 \sum_{i:\gamma_i>0}\gamma_i|Z^{\gamma-e_i}|\le d.
\]
Consequently,
\[
 |\nabla F_t|\ge d^{(4-d)/2}-d|t|_1
 \ge\frac12\quad(d=4),
 \qquad
 |\nabla F_t|>\frac7{80}\quad(d=5).
\]
The bounds remain positive on a slightly larger open coefficient
ball $U_d\supset B_d$. Thus $dF_t$ never vanishes on
$F_t^{-1}(0)\setminus\{0\}$, and every $X_t$ is nonsingular.

Consider the projection
\[
 \pi:\{(t,[Z])\in U_d\times\PP^{d-1}:F_t(Z)=0\}
          \longrightarrow U_d,\qquad (t,[Z])\longmapsto t.
\]
The lower bound for the derivatives makes $\pi$ a submersion.
The projection is proper because $\PP^{d-1}$ is compact.
A smooth horizontal distribution therefore lifts each straight
path from $0$ to $t$ for its full length. Transport along
the path and its reverse gives a diffeomorphism $X_0\simeq X_t$.
Since the Fermat hypersurface $X_0$ is connected, so is $X_t$.
\end{proof}

For $X_t\subset\PP^{d-1}$, the adjunction formula gives
\[
 K_{X_t}
 \simeq
 \bigl(K_{\PP^{d-1}}\otimes\mathcal O_{\PP^{d-1}}(d)\bigr)|_{X_t}
 \simeq\mathcal O_{X_t}.
\]
For the adjunction formula and the Poincar\'e residue map, see, e.g.,
Griffiths and Harris
\cite[Chapter~1, Section~1, Adjunction Formula~II, p.~147]{GH}.
The residue map gives the nowhere vanishing holomorphic $n$-form
\[
 \Omega_t=\operatorname{Res}_{X_t}
 \frac{\sum_{a=0}^{d-1}(-1)^a Z_a\,
       dZ_0\wedge\cdots\wedge\widehat{dZ_a}\wedge\cdots\wedge dZ_{d-1}}
      {F_t}.
\]
For a K\"ahler form $\alpha\in[\omega_{0,t}]$, define
\begin{equation}
 dV_{\Omega_t}=(\sqrt{-1})^{n^2}\Omega_t\wedge\overline\Omega_t,\qquad
 R_\alpha=\frac{\alpha^n/n!}{dV_{\Omega_t}},\qquad
 C_t=\frac{\int_{X_t}\omega_{0,t}^n/n!}{\int_{X_t}dV_{\Omega_t}}.
 \label{eq:CY-density}
\end{equation}
The Ricci-flat form satisfies $\beta_t^n/n!=C_t\,dV_{\Omega_t}$,
so the normalized measure $d\nu_t$ defined above is
\[
 d\nu_t=\frac{dV_{\Omega_t}}{\int_{X_t}dV_{\Omega_t}}.
\]
In particular,
\begin{equation}
 q=\frac{\beta_t^n}{\alpha^n}=\frac{C_t}{R_\alpha},\qquad
 C_t=\int_{X_t}R_\alpha\,d\nu_t.
 \label{eq:normalized-density}
\end{equation}
The ratio $\sup R_\alpha/\inf R_\alpha$ does not depend on
rescaling $\Omega_t$. The Ricci form is
\begin{equation}
 \operatorname{Ric}(\alpha)
 =-\sqrt{-1}\partial_{X_t}\bar\partial_{X_t}\log R_\alpha.
 \label{eq:ricci-residue-density}
\end{equation}

\subsection{Local expressions and polynomial inequalities}
By a polynomial $P(t,\Pi)$ with rational coefficients we mean
a polynomial in the real variables
$\operatorname{Re}t_\gamma$, $\operatorname{Im}t_\gamma$,
$\operatorname{Re}\Pi_{ab}$ and $\operatorname{Im}\Pi_{ab}$
whose coefficients lie in $\Q$. Restricting the polynomial to
$\{(t,[Z]):t\in B_d,\ [Z]\in X_t\}$ gives a smooth real-valued
function. Put
\[
 \alpha_{P,t}=\omega_{0,t}
       +\sqrt{-1}\partial_{X_t}\bar\partial_{X_t}P(t,\Pi),
 \qquad R_P=R_{\alpha_{P,t}}.
\]

In the affine chart $Z_a=1$, write $w_i=Z_i/Z_a$ for $i\ne a$
and $f_a(t,w)=F_t|_{Z_a=1}$. At a point where
$\partial f_a/\partial w_b\ne0$, the holomorphic implicit
function theorem expresses $w_b$ as a holomorphic function of
the $n$ coordinates $(w_j)_{j\notin\{a,b\}}$ on $X_t$.
Write $(f_a)_j=\partial f_a/\partial w_j$. The coordinate
tangent vectors and derivatives are
\begin{equation}
 T_{ab,j}=e_j-\frac{(f_a)_j}{(f_a)_b}e_b,
 \qquad
 D_j=\partial_{w_j}-\frac{(f_a)_j}{(f_a)_b}\partial_{w_b}
 \quad(j\notin\{a,b\}),
 \label{eq:coordinate-tangent-derivatives}
\end{equation}
where $e_j$ are the coordinate vectors of the affine chart.
With $S=1+\sum_{i\ne a}|w_i|^2$, set
\[
 G_a=\bigl(\partial_{\bar w_i}\partial_{w_j}
                       (\log S+P(t,\Pi))\bigr)_{i,j\ne a},
 \qquad g_{ab}=T_{ab}^*G_aT_{ab}.
\]
In the coordinates given by the implicit function theorem,
the residue form is, up to a sign,
$\Omega_t=\bigwedge_{j\notin\{a,b\}}dw_j/(f_a)_b$. Hence
\begin{equation}
 R_P(t,w)=\det(g_{ab})\,|(f_a)_b|^2.
 \label{eq:local-density}
\end{equation}
The entries of $g_{ab}$ and $R_P$ are rational functions of
$\operatorname{Re}w_i$, $\operatorname{Im}w_i$,
$\operatorname{Re}t_\gamma$ and $\operatorname{Im}t_\gamma$.
Their denominators are products of powers of $S$ and
$|(f_a)_b|^2$, both positive on the coordinate domain.

For background on real algebraic geometry, see, for example,
Bochnak, Coste and Roy \cite{BCR}. We use the following decision
theorem.
\begin{lemma}[{Basu--Pollack--Roy \cite[Chapter~14, Section~14.1]{BPR}}]
\label{lem:real-elimination}
Every first-order sentence over the real numbers formed from
polynomial equalities and inequalities with rational coefficients
is decidable by a finite algorithm.
\end{lemma}

For a K\"ahler form
$\alpha=\sqrt{-1}h_{i\bar j}\,dz_i\wedge d\bar z_j$, we use
the curvature tensor
\begin{equation}
 R^\alpha_{i\bar j k\bar l}
 =-\partial_i\partial_{\bar j}h_{k\bar l}
   +h^{p\bar q}(\partial_i h_{k\bar q})
                    (\partial_{\bar j}h_{p\bar l}).
 \label{eq:curvature-convention}
\end{equation}
The notation $\operatorname{Bisec}(\alpha)\le K$ means
\[
 R^\alpha_{i\bar j k\bar l}
        \xi^i\overline{\xi^j}\zeta^k\overline{\zeta^l}
 \le K|\xi|_\alpha^2|\zeta|_\alpha^2
 \quad\text{for all }\xi,\zeta\in T^{1,0}X.
\]
Here $R^\alpha$ denotes curvature, whereas $R_\alpha$ in
\eqref{eq:CY-density} denotes a scalar volume density.

\begin{lemma}\label{lem:QE}
Let $P(t,\Pi)$ be a polynomial with rational coefficients and
let $0<\ell\le1\le M$, $K\ge0$ and $\kappa>1$ be rational.
One can decide in finitely many operations on rational numbers
whether
\begin{align}
 &\ell\omega_{0,t}\le\alpha_{P,t}\le M\omega_{0,t},
       \qquad\operatorname{Bisec}(\alpha_{P,t})\le K,
       \label{eq:parametric-geometric-tests}\\
 &R_P(t,x)\le\kappa R_P(t,y)\qquad(x,y\in X_t)
       \label{eq:parametric-density-test}
\end{align}
hold for every $t\in B_d$. For an integer $s\ge1$ and a
rational number $\eta>0$, one can also decide whether
\begin{equation}
 \sup_{t\in B_d}\sup_{X_t}
 |\nabla_{0,t}^{\,k}\log R_P|_{g_{0,t}}<\eta
 \quad(1\le k\le s).
 \label{eq:log-density-derivative-tests}
\end{equation}
\end{lemma}
\begin{proof}
First express $B_d$ by polynomial conditions. Introduce real
variables $s_\gamma$ satisfying
\[
 s_\gamma\ge0,\qquad
 s_\gamma^2=(\operatorname{Re}t_\gamma)^2
                +(\operatorname{Im}t_\gamma)^2,\qquad
 \sum_\gamma s_\gamma\le\rho_d,
 \quad\rho_4=\tfrac18,\quad\rho_5=\tfrac1{16}.
\]
On each coordinate domain impose
$\operatorname{Re}f_a=\operatorname{Im}f_a=0$ and
$|(f_a)_b|^2>0$. These domains cover $X_t$. To see this,
suppose that all affine derivatives vanished at a point with $Z_a\ne0$.
Then
\[
 \sum_i Z_i\frac{\partial F_t}{\partial Z_i}=dF_t=0
\]
would also give $\partial F_t/\partial Z_a=0$, contradicting
Lemma~\ref{lem:family-smoothness}.

Let $g_{0,ab}$ be the matrix of $\omega_{0,t}$ in these
coordinates. The metric inequalities hold exactly when all
principal minors of $g_{ab}-\ell g_{0,ab}$ and
$Mg_{0,ab}-g_{ab}$ are nonnegative. The decision procedure
starts with these conditions. If they hold, then $\alpha_{P,t}>0$,
so the subsequent curvature and logarithmic derivatives are defined.

To test the curvature bound, apply $D_j$ from
\eqref{eq:coordinate-tangent-derivatives}
and its complex conjugate to the entries of $g_{ab}$.
Formula \eqref{eq:curvature-convention} then expresses
curvature by rational functions in the same real variables.
The inequality $\operatorname{Bisec}(\alpha_{P,t})\le K$
becomes a rational inequality quantified over the real and
imaginary components of two tangent vectors.

To test \eqref{eq:parametric-density-test}, use
\eqref{eq:local-density} on each pair of coordinate domains
with the same parameter $t$.

The inverse background metric and its Christoffel symbols
also have rational coefficients. Starting with
$d\log R_P=R_P^{-1}dR_P$, repeated use of the covariant
derivative formula gives rational coefficients for
$\nabla_{0,t}^{\,k}\log R_P$ for every $k\ge1$. The squared
norm of each tensor, computed using the inverse real metric,
is therefore rational. Since $\eta>0$, the required
pointwise bound can be written as
\[
 |\nabla_{0,t}^{\,k}\log R_P|_{g_{0,t}}^2<\eta^2.
\]

For each order of differentiation, the norm is continuous on the compact set
$\{(t,[Z]):t\in B_d,\ [Z]\in X_t\}$; hence the pointwise
strict inequalities are equivalent to the strict supremum
bounds in \eqref{eq:log-density-derivative-tests}.

Finally, replace each rational inequality $a/b\ge0$ by
$ab\ge0$, and each $a/b>0$ by $ab>0$. Each replacement
multiplies the inequality by $b^2>0$. The coordinate conditions
and the established positivity of the metrics and $R_P$ ensure
that the denominators do not vanish. Thus all the required
conditions are finite statements of polynomial equalities and
inequalities with rational coefficients.
Lemma~\ref{lem:real-elimination} decides them.
\end{proof}

\section{Explicit metrics on the Fermat hypersurfaces}
\label{sec:explicit-potential}
We give the rational-coefficient polynomials in
Theorem~\ref{thm:explicit-main} and prove their integral bounds.
The computer-assisted argument in Lemma~\ref{lem:certified-density} gives global density bounds;
the comparison with the Ricci-flat forms is analytic.

\subsection{The explicit potentials}
\label{subsec:recovered-potentials}
For $d=4,5$, put
\[
 x_i=\frac{|Z_i|^2}{\sum_k|Z_k|^2},\qquad
 s_j=\sum_{i=0}^{d-1}x_i^j.
\]
The finite group
\[
 G_d=\bigl((\mu_d)^d/\mu_d\bigr)\rtimes\mathfrak S_d
\]
acts by permuting the coordinates and multiplying them by $d$th
roots of unity. The quotient identifies transformations that differ
by common scalar multiplication. This action preserves $X_d$ and
$\omega_0$, and the functions $s_j$ are invariant under it.

The polynomials $P_{Q_{12}}$ and $P_{T_{12}}$ are defined in
\eqref{eq:recovered-Q12-T12} and Appendix~\ref{app:coefficients}.
We write $P_4=P_{Q_{12}}$, $P_5=P_{T_{12}}$, and
\begin{equation}
 \alpha_A=\omega_0+\sqrt{-1}\partial\bar\partial P_A\big|_{X_d}.
 \label{eq:recovered-forms}
\end{equation}
Each $P_A$ is a smooth real function on projective space.
The form $\alpha_A$ has local potential
$\log\sum_i|Z_i|^2+P_A$.
Since $s_j=\sum_i\Pi_{ii}^j$, both polynomials have degree
at most twelve in $\Pi$.

For comparison, the lower-degree polynomial for the quartic is
\begin{equation}
\begin{aligned}
10^7P_{Q_8}={}&-71074s_2+31334581s_3-37805554s_4
                 +52868708s_2^2\\
 &-301613657s_2s_3+321972842s_3^2+371696981s_2s_4\\
 &-85441278s_2^3-789158771s_3s_4+273188204s_2^2s_3\\
 &+496612443s_4^2-26380666s_2s_3^2
                 -320501478s_2^2s_4+32661183s_2^4.
\end{aligned}
\label{eq:recovered-Q8}
\end{equation}
It has fourteen terms and degree at most eight in $\Pi$.
All three resulting real $(1,1)$ forms belong to $[\omega_0]$.

\subsection{Global volume-density bounds}
We use $R_\alpha=(\alpha^n/n!)/dV_\Omega$ also when $\alpha$
has not yet been shown to be positive. Formula
\eqref{eq:local-density} computes this signed density.

\begin{lemma}[Global density bounds; computer-assisted]
\label{lem:certified-density}
The forms $\alpha_{Q_8}$, $\alpha_{Q_{12}}$ and $\alpha_{T_{12}}$
satisfy the following density bounds on the indicated Fermat
hypersurfaces:
\begin{equation}
\begin{array}{c|c|cc}
 A&d&L&U\\ \hline
 Q_8&4&333/200&337/200\\
 Q_{12}&4&1667/1000&837/500\\
 T_{12}&5&537/1000&139/250.
\end{array}
\label{eq:concrete-density-bounds}
\end{equation}
Here each row means $L\le R_{\alpha_A}\le U$ on $X_d$.
\end{lemma}
\begin{proof}
The proof reduces each density bound to finitely many
inequalities between integers. We give the polynomial identities
and the subdivision procedure used to verify those inequalities.

For this calculation, use the unnormalized variables
$x_i=|Z_i|^2$ and $S=\sum_i x_i$. Homogenize each potential as
$P=Q/(DS^m)$, where $Q$ is an integer homogeneous polynomial
of degree $m$ and $D>0$ is an integer. Differentiating
$K=\log S+P$ gives
\begin{equation}
 w_{ij}=\frac{DS^m-S^2Q_{ij}+mS(Q_i+Q_j)-m(m+1)Q}
                  {DS^{m+2}},\qquad i<j,
 \label{eq:certificate-weights}
\end{equation}
where $Q_i=\partial Q/\partial x_i$ and
$Q_{ij}=\partial^2Q/\partial x_i\partial x_j$.
Here $w_{ij}$ is the coefficient of the square
$|Z_iV_j-Z_jV_i|^2$ in the Hermitian Hessian $G$ of $K$:
\[
 G[V,\overline V]=\sum_{i<j}w_{ij}|Z_iV_j-Z_jV_i|^2.
\]
After the coordinate factors are removed, the off-diagonal
coefficients are $-w_{ij}$. Homogeneity then determines the
diagonal coefficients through $GZ=0$.

Choose representatives with $S=1$. Let $\mathcal F_2$ be the
set of graphs on $\{0,\ldots,d-1\}$ which have two connected
components and no cycles. Each $T\in\mathcal F_2$ has $d-2$
edges; write $E(T)$ for its edge set, $A,B$ for its components,
and $\deg_T(i)$ for the number of edges incident to $i$. Then
\begin{equation}
 R_\alpha=d^2\sum_{T\in\mathcal F_2}
 \left(\prod_{ij\in E(T)}w_{ij}\right)
 \left(\prod_i x_i^{\deg_T(i)-1}\right)
 \left|\sum_{i\in A}Z_i^d\right|^2.
 \label{eq:certificate-forest}
\end{equation}
To prove this identity, put
\[
 f=(dZ_i^{d-1})_i,\qquad
 H_Z=\{V:Z^*V=0,\ fV=0\}.
\]
Let $U$ be a $d\times n$ matrix whose columns form a Euclidean
orthonormal basis of $H_Z$. Formula \eqref{eq:local-density},
with the residue form defined in Section~\ref{sec:preliminaries},
reads
\[
 R_\alpha=\|f\|^2\det(U^*GU).
\]
For a set $T$ of $n$ edges, let $E_T$ be the matrix with rows
$V\mapsto Z_iV_j-Z_jV_i$, $ij\in T$. Expanding the
determinant of the sum of these rank-one matrices gives
\[
 \det(U^*GU)
 =\sum_{|T|=n}\left(\prod_{ij\in T}w_{ij}\right)
       |\det(E_TU)|^2.
\]
This determinant identity does not require $w_{ij}>0$.
If the chosen edges contain a cycle, the corresponding rows are
linearly dependent. Thus only $T\in\mathcal F_2$ contribute.

For such a forest, put $s_A=\sum_{i\in A}x_i$ and
$s_B=\sum_{i\in B}x_i$. Initially assume all $x_i>0$.
The matrix $E_TE_T^*$ consists of the pairwise Hermitian
inner products of the rows of $E_T$. Its determinant is
\begin{equation}
 \det(E_TE_T^*)=s_As_B\prod_i x_i^{\deg_T(i)-1}.
 \label{eq:forest-edge-determinant}
\end{equation}
Delete the column of a vertex
$r\in A$ from the edge matrix of the tree $A$. Successively
removing leaves shows that the squared modulus of this minor is
$x_r\prod_{i\in A}x_i^{\deg_T(i)-1}$. Summing over $r\in A$
and repeating the calculation for $B$ gives the formula.

The kernel of $E_T$ is spanned by $Z|_A$ and $Z|_B$.
Write $f^*=(d\overline Z_i^{d-1})_i$ as a column vector, and
let $\pi_T$ be the Euclidean orthogonal projection onto
$\ker E_T$. Since $\sum_iZ_i^d=0$ and $s_A+s_B=1$,
\[
 \pi_T f^*\perp Z,\qquad
 \|\pi_T f^*\|^2
 =d^2\left(\frac1{s_A}+\frac1{s_B}\right)
      \left|\sum_{i\in A}Z_i^d\right|^2
 =\frac{d^2}{s_As_B}\left|\sum_{i\in A}Z_i^d\right|^2.
\]
Computing the determinant in orthonormal bases of
$\ker E_T$ and its orthogonal complement gives
\[
 \|f\|^2|\det(E_TU)|^2
 =\det(E_TE_T^*)\|\pi_T f^*\|^2.
\]
Combining this equality with \eqref{eq:forest-edge-determinant}
proves \eqref{eq:certificate-forest}. If $i$ is an isolated
vertex, choose $A=\{i\}$. The factor $x_i^{-1}$ is then cancelled
by $|Z_i^d|^2=x_i^d$. All remaining exponents of $x_i$ are
nonnegative, so continuity extends the formula to points with
zero coordinates.

For the expansion below, choose $A$ to be the component with
fewer vertices; if both components have two vertices, choose
the one containing $0$. Expanding \eqref{eq:certificate-forest} gives
\[
 R_\alpha=\frac{D_0+2\sum_{i<j}A_{ij}r_{ij}}{H},
 \qquad r_{ij}=\operatorname{Re}(Z_i^d\overline{Z_j^d}),
 \qquad H>0.
\]
Compute $D_0$ and $A_{ij}$ by substituting
\eqref{eq:certificate-weights} into \eqref{eq:certificate-forest},
multiplying by $D^nS^{n(m+2)}$, and expanding
$|\sum_{i\in A}Z_i^d|^2$. These operations give a polynomial
numerator; on $S=1$ the denominator is $H=D^n$.

Next multiply $\sum_iZ_i^d=0$ by $\overline{Z_i^d}$
and take real parts to obtain $\sum_{j\ne i}r_{ij}=-x_i^d$.
Consequently, with
\[
 a_0=A_{01}+A_{02}-A_{12},\qquad
 a_j=2A_{0j}-a_0\ (j>0),\qquad
 E_{ij}=2A_{ij}-a_i-a_j,
\]
the numerator equals
\begin{equation}
 D_0-\sum_i a_ix_i^d+\sum_{i<j}E_{ij}r_{ij},
 \qquad E_{0j}=E_{12}=0.
 \label{eq:certificate-phases}
\end{equation}
For $d=4$ use $|r_{ij}|\le x_i^2x_j^2$; for $d=5$ use
$|r_{ij}|\le x_i^2x_j^2(x_i+x_j)/2$.
The terms still involving $r_{ij}$ have indices
$(i,j)\in\{(1,3),(2,3)\}$ when $d=4$, and
$(i,j)\in\{(1,3),(1,4),(2,3),(2,4),(3,4)\}$ when $d=5$.

By symmetry it suffices to consider
\[
 x_0\ge x_1\ge\cdots\ge x_{d-1}\ge0,
 \qquad \sum_i x_i=1.
\]
The equation $Z_0^d=-\sum_{i>0}Z_i^d$ and the triangle
inequality imply
\begin{equation}
 x_0^{d/2}\le\sum_{i>0}x_i^{d/2}\le(1-x_0)^{d/2};
 \qquad x_0\le\tfrac12.
 \label{eq:fermat-magnitude-constraint}
\end{equation}
For an explicit covering, set
\[
 y_k=k(x_{k-1}-x_k)\ (1\le k<d),\qquad y_d=dx_{d-1}.
\]
Then $y_k\ge0$, $\sum_k y_k=1$, and
$x_i=\sum_{k=i+1}^d y_k/k$. The cap condition becomes
\[
 y_1\le\sum_{k=3}^d b_ky_k,\qquad b_k=(k-2)/k.
\]
Partition this region by
\[
 \sum_{k=3}^{r-1}b_ky_k\le y_1\le
 \sum_{k=3}^r b_ky_k,\qquad 3\le r\le d.
\]
The corresponding cone has independent rays
\[
 e_2,\qquad b_ke_1+e_k\ (3\le k\le r),\qquad
e_k\ (r\le k\le d).
\]
Here $e_k$ are the coordinate vectors. Normalize each ray so
that its coordinates sum to one, then apply the inverse formula
for $x$. The resulting $d-2$ rational simplices have disjoint
interiors and cover the capped ordered domain. Order their vertices
as in the displayed list of rays, and order the simplices by
increasing $r$.

The refined quartic calculation also uses
the following exact relation. Put $r_i=x_i^2$ and
\[
 t_{ij}=r_ir_j-\operatorname{Re}(Z_i^4\overline{Z_j^4})
 \quad(1\le i<j\le3).
\]
Then $t_{ij}\ge0$ and
\[
 t_{12}+t_{13}+t_{23}
 =\Delta:=\frac{(x_1^2+x_2^2+x_3^2)^2-x_0^4}{2}.
\]
By \eqref{eq:certificate-phases}, the density numerator is
$H R_\alpha=B-E_{13}t_{13}-E_{23}t_{23}$, where
\[
 B=D_0-\sum_i a_ix_i^4+E_{13}x_1^2x_3^2+E_{23}x_2^2x_3^2.
\]
Thus $\Delta\ge0$ at every point of $X_4$. Define
\begin{equation}
 B_0=B,\qquad B_1=B-E_{13}\Delta,\qquad
 B_2=B-E_{23}\Delta.
 \label{eq:quartic-density-polynomials}
\end{equation}
If $\Delta>0$, then
\[
 H R_\alpha=\frac{t_{12}}{\Delta}B_0+\frac{t_{13}}{\Delta}B_1
       +\frac{t_{23}}{\Delta}B_2,\qquad
 \frac{t_{12}+t_{13}+t_{23}}{\Delta}=1.
\]
Hence $\min_{0\le k\le2}B_k\le H R_\alpha\le\max_{0\le k\le2}B_k$;
when $\Delta=0$, $H R_\alpha=B_0$.

Here is the finite arithmetic test. Put $M=2^{40}$ and
$N=n(m+2)$, so $N=20,28,42$ for $Q_8,Q_{12},T_{12}$,
respectively. Set
\[
 F_0=\frac{D_0-\sum_i a_ix_i^d}{H},\qquad
 F_{ij}=\frac{E_{ij}b_{ij}}{H},\qquad
 b_{ij}=\begin{cases}
 x_i^2x_j^2,&d=4,\\
 x_i^2x_j^2(x_i+x_j)/2,&d=5.
 \end{cases}
\]
For $b_{ij}>0$, put $\xi_{ij}=r_{ij}/b_{ij}$, and put
$\xi_{ij}=0$ when $b_{ij}=0$. Then
$R_\alpha=F_0+\sum_{i<j}F_{ij}\xi_{ij}$ and
$|\xi_{ij}|\le1$. For $Q_{12}$ also use
$B_0/H,B_1/H,B_2/H$ from \eqref{eq:quartic-density-polynomials}.
All these polynomials are homogeneous of degree $N$.

On a simplex with ordered vertices $v_0,\ldots,v_{d-1}$,
substitute $x=\sum_a\theta_av_a$. If
\[
 F(x)=\sum_{|e|=N}a_e\theta^e
     =\sum_{|e|=N}c_e\binom Ne\theta^e,
 \qquad c_e=\frac{a_e}{\binom Ne},
\]
the functions $\binom Ne\theta^e$, $|e|=N$, are its
Bernstein basis functions. Here $\theta_a\ge0$ and
$\sum_a\theta_a=1$, so the multinomial formula gives
$\sum_{|e|=N}\binom Ne\theta^e=1$ and therefore
$\min_ec_e\le F\le\max_ec_e$.
For each initial simplex compute the rational coefficients
exactly and put
\[
 l_e=\lfloor Mc_e\rfloor,\qquad u_e=\lceil Mc_e\rceil.
\]
Write $l_{0,e},u_{0,e}$ for the bounds of $F_0$.
List the nonzero polynomials $F_{ij}$ in lexicographic order
of $(i,j)$, and denote their coefficient bounds by
$l_{p,e},u_{p,e}$. Set $L_M=\lceil ML\rceil$ and
$U_M=\lfloor MU\rfloor$.
The density test accepts a simplex if, for every $e$,
\begin{equation}
 l_{0,e}-\sum_p\max\{|l_{p,e}|,|u_{p,e}|\}\ge L_M,
 \qquad
 u_{0,e}+\sum_p\max\{|l_{p,e}|,|u_{p,e}|\}\le U_M.
 \label{eq:certificate-acceptance}
\end{equation}
For $Q_{12}$, one may instead verify that every lower and
upper coefficient bound of $B_k/H$, $k=0,1,2$, lies in
$[L_M,U_M]$. Either test proves $L\le R_\alpha\le U$
throughout that simplex.

For $Q_8$ and $T_{12}$, simplify the phase bounds before
\eqref{eq:certificate-acceptance}. If $l_{p,e}\ge0$ for every $e$,
replace $(l_{0,e},u_{0,e})$ by
$(l_{0,e}-u_{p,e},u_{0,e}+u_{p,e})$ and remove $p$.
If $u_{p,e}\le0$ for every $e$, replace them by
$(l_{0,e}+l_{p,e},u_{0,e}-l_{p,e})$ and remove $p$.
Process $p$ in decreasing order. These operations bound the
contribution of $F_{ij}\xi_{ij}$ for every $|\xi_{ij}|\le1$.
Subsequent subdivisions use the resulting bounds.
No such simplification is made for $Q_{12}$.

Before testing the density, use the following criteria to discard
simplices containing no point satisfying
\eqref{eq:fermat-magnitude-constraint}. For $d=4,5$, this
constraint implies
$x_0^2\le(\sum_{i>0}x_i^{d/2})^{4/d}\le\sum_{i>0}x_i^2$.
The first exclusion test is therefore
\begin{equation}
 2v_{a0}v_{b0}-\sum_jv_{aj}v_{bj}>0
 \quad\text{for every }0\le a\le b\le d-1.
 \label{eq:certificate-quadratic-exclusion}
\end{equation}
These are the quadratic Bernstein coefficients of
$x_0^2-\sum_{j>0}x_j^2$. For $d=5$, put
\[
 a=\frac1d\sum_bv_{b0},\qquad
 \ell(z)=\frac{\lfloor M\sqrt z\rfloor}{M},\qquad
 u(z)=\ell(z)+\frac1M.
\]
The square-root bounds are determined by integer square
comparisons. The second exclusion test requires, at every vertex,
\[
 5v_{b0}-3a\ge0,\qquad
 \frac12a\ell(a)(5v_{b0}-3a)
       >\sum_{j>0}v_{bj}^2u(v_{bj}).
\]
Indeed, convexity gives the tangent-line bound
$x_0^{5/2}\ge a^{3/2}(5x_0-3a)/2$, and interpolation of
$\sum_{j>0}x_j^{5/2}$ at the vertices gives an affine upper
bound. The test therefore gives
$x_0^{5/2}>\sum_{j>0}x_j^{5/2}$ throughout the simplex,
contrary to \eqref{eq:fermat-magnitude-constraint}.

If a simplex is neither excluded nor accepted by the density
test, bisect its longest Euclidean edge, resolving ties
lexicographically. For an edge $(i,j)$, let $v_*=(v_i+v_j)/2$.
The first subsimplex replaces $v_j$ by $v_*$, and the second
replaces $v_i$ by $v_*$.
To obtain their coefficient bounds, fix all exponents except
$e_i,e_j$ and write $e_i+e_j=k$. Index these coefficients by
$e_j=q$, $0\le q\le k$. Starting from the existing bounds
$l_q^{(0)},u_q^{(0)}$, use
\[
 l_q^{(h+1)}=\left\lfloor\frac{l_q^{(h)}+l_{q+1}^{(h)}}2\right\rfloor,
 \qquad
 u_q^{(h+1)}=\left\lceil\frac{u_q^{(h)}+u_{q+1}^{(h)}}2\right\rceil.
\]
The first subsimplex has coefficient bounds with index $q$
$(l_0^{(q)},u_0^{(q)})$; the second has
$(l_q^{(k-q)},u_q^{(k-q)})$.
These averages come from substituting the midpoint in the
polynomial expansion. Rounding down for $l$ and up for $u$
preserves the lower and upper bounds. Process the first
subsimplex before the second.

The calculation succeeds only if every terminal simplex is
excluded or passes the density test. It fails if a simplex
remains undecided after $52$ bisections or if more than
$500000$ simplices are examined.

For the displayed coefficients, the seven initial-simplex calculations finish
with the following counts:
\[
\begin{array}{c|rrrr}
 A&\text{root simplices}&\text{nodes}&\text{terminal cells}
       &\text{maximum depth}\\ \hline
 Q_8&2&3824&1913&21\\
 Q_{12}&2&5686&2844&22\\
 T_{12}&3&14295&7149&27.
\end{array}
\]
The counts record exact integer calculations. Every terminal
simplex either satisfies $L\le R_\alpha\le U$ or contains
no point of $X_d$. Since these simplices cover the ordered
domain, they prove \eqref{eq:concrete-density-bounds}.
\end{proof}

\subsection{The analytic comparison}
\begin{lemma}[Positivity from the volume form]
\label{lem:positive-volume}
Let $(X,\omega_0)$ be a compact connected K\"ahler manifold, and let $P$ be
a smooth real function. If
$(\omega_0+\sqrt{-1}\partial\bar\partial P)^n$ never vanishes, then
$\omega_0+\sqrt{-1}\partial\bar\partial P$ is positive.
\end{lemma}
\begin{proof}
At a minimum of $P$, its complex Hessian is nonnegative, so the
form is positive there. The eigenvalues relative to $\omega_0$
are continuous and never zero. Connectedness therefore forces
them to remain positive throughout $X$.
\end{proof}

\begin{lemma}[Integral tensor comparison]\label{lem:volume-contrast}
Let $\alpha,\beta$ be cohomologous K\"ahler forms on a compact
complex $n$-fold, and put
\[
 \rho=\frac{\sup_X(\beta^n/\alpha^n)}
             {\inf_X(\beta^n/\alpha^n)}.
\]
Then
\begin{equation}
 \max_{\gamma\in\{\alpha,\beta\}}
 \int_X|\alpha-\beta|_{\gamma,\HS}\,d\mu_\gamma
 \le\frac{\sqrt n}{2}\log\rho.
 \label{eq:volume-contrast-L1}
\end{equation}
\end{lemma}
\begin{proof}
Let $\lambda_i$ be the positive eigenvalues of
$\alpha^{-1}\beta$ and put $q=\prod_i\lambda_i$.
Cohomology and Stokes' theorem give
\[
 \int_Xq\,d\mu_\alpha=1,\qquad
 \int_X\sum_i\lambda_i\,d\mu_\alpha
 =\frac{n\int_X\beta\wedge\alpha^{n-1}}
         {\int_X\alpha^n}=n.
\]
Set $Y=\log q$, $\ell=\log\rho$ and
\[
 \psi(t)=\log\int_Xe^{tY}\,d\mu_\alpha,
 \qquad d\mu_t=e^{tY-\psi(t)}d\mu_\alpha.
\]
Differentiation gives
\[
 \psi''(t)=\int_X\left(Y-\int_XY\,d\mu_t\right)^2d\mu_t
 \le\frac{\ell^2}{4}.
\]
Indeed, if $a$ is the midpoint of the range of $Y$, then
$\operatorname{Var}_{\mu_t}(Y)\le\int(Y-a)^2d\mu_t\le\ell^2/4$.
Since $\psi(0)=\psi(1)=0$, Taylor's formula for $\psi$ gives
\[
 H:=-\int_X\log q\,d\mu_\alpha
   =\int_0^1(1-t)\psi''(t)\,dt\le\frac{\ell^2}{8}.
\]

On $X\times\{1,\ldots,n\}$ define
\[
 d\sigma(x,i)=\frac1n\,d\mu_\alpha(x),\qquad
 h(x,i)=\lambda_i(x).
\]
Then $\int d\sigma=\int h\,d\sigma=1$ and
$-\int\log h\,d\sigma=H/n$. We claim that
\begin{equation}
 \int|h-1|\,d\sigma
 \le\left(-2\int\log h\,d\sigma\right)^{1/2}.
 \label{eq:scalar-L1-log}
\end{equation}
To prove \eqref{eq:scalar-L1-log}, take $A=\{h\le1\}$,
$a=\sigma(A)$, and $b=\int_Ah\,d\sigma$.
Then $\int|h-1|\,d\sigma=2(a-b)$, and Jensen's inequality
on $A$ and its complement gives
\[
 -\int\log h\,d\sigma
 \ge a\log\frac ab+(1-a)\log\frac{1-a}{1-b}
 \ge2(a-b)^2.
\]
For fixed $b\in(0,1)$, put
\[
 F_b(a)=a\log\frac ab+(1-a)\log\frac{1-a}{1-b}.
\]
Then $F_b(b)=F_b'(b)=0$ and $F_b''(a)=1/(a(1-a))\ge4$,
so $F_b(a)\ge2(a-b)^2$. The endpoint cases follow by continuity.
This proves \eqref{eq:scalar-L1-log}. Consequently,
\[
 \int_X|\alpha-\beta|_{\alpha,\HS}\,d\mu_\alpha
 \le n\int|h-1|\,d\sigma
 \le\sqrt{2nH}\le\frac{\sqrt n}{2}\ell.
\]
Interchanging $\alpha$ and $\beta$ replaces $q$ by $q^{-1}$
and leaves its logarithmic range unchanged. The cohomological
identities also hold for $\mu_\beta$, so the preceding argument
gives the other direction of \eqref{eq:volume-contrast-L1}.
\end{proof}

\begin{proof}[Proof of Theorem~\ref{thm:explicit-main}]
The density lower bounds in Lemma~\ref{lem:certified-density}
and Lemma~\ref{lem:positive-volume} show that
$\widehat\omega_4=\alpha_{Q_{12}}$ and
$\widehat\omega_5=\alpha_{T_{12}}$ are K\"ahler.
The exact Ricci-flat volume form is a constant multiple of
the residue volume form, so the volume ratio in
Lemma~\ref{lem:volume-contrast} satisfies
\[
 \rho\le\frac{1674}{1667}\quad(d=4),
 \qquad \rho\le\frac{556}{537}\quad(d=5).
\]
Since $\log(1+x)<x$ for $x>0$, the respective integral bounds
in \eqref{eq:volume-contrast-L1} are smaller than
\[
 \frac{7\sqrt2}{3334}<\frac3{1000},\qquad
 \frac{19\sqrt3}{1074}<\frac{31}{1000}.
\]
This proves \eqref{eq:explicit-main-L1}.
\end{proof}

The same argument applies to $Q_8$, with
$\rho\le337/333$, so this lower-degree formula is also
a K\"ahler approximation with an evaluated integral bound.
The normalization of the exact residue volume is not needed
in these comparisons.

\section{Uniform estimates for the projective background}\label{sec:projective-background}\label{app:background}
We determine Sobolev and Poincar\'e constants for
$(X_t,\omega_{0,t})$ that are uniform for $t\in B_d$.
The constants are obtained from rational inequalities for
the defining equations of the hypersurfaces.

\subsection{A Sobolev constant}
For submanifold Sobolev inequalities, see, for example,
Michael and Simon \cite{MichaelSimon} and Hoffman and Spruck
\cite{HoffmanSpruck}. We use the following explicit constant.
\begin{lemma}[{Brendle \cite[Theorem~1]{Brendle}}]
\label{lem:brendle}
Let $M^m$ be a compact submanifold without boundary in
$\mathbb R^{m+k}$, with $k\ge2$. Write $H$ for the trace of its
second fundamental form, and let $B^q$ denote the unit Euclidean
ball in dimension $q$. Every smooth positive function $a$ satisfies
\[
 m\left(\frac{(m+k)|B^{m+k}|}{k|B^k|}\right)^{1/m}
 \|a\|_{m/(m-1),dV}
 \le\int_M\sqrt{|\nabla a|^2+|H|^2a^2}\,dV.
\]
\end{lemma}

\begin{proposition}[Projective Sobolev bound]
\label{prop:projective-sobolev}
Let $X^n\subset\mathbb P^{n+1}$ be a smooth quartic surface
or quintic threefold, and let $\omega_0$ be its restricted
Fubini--Study form. Every smooth real-valued function $w$ satisfies
\begin{equation}
 \|w\|_{2n/(n-1),\mu_0}^2
 \le102\bigl(\|\partial w\|_{2,\mu_0}^2+\|w\|_{2,\mu_0}^2\bigr),
 \qquad d\mu_0=\frac{\omega_0^n}{\int_X\omega_0^n}.
 \label{eq:appendix-projective-S}
\end{equation}
\end{proposition}
\begin{proof}
Put $N=d=n+2$, $m=2n$, $D=N^2-1$, and $k=D-m$.
The map
\[
 P([Z])=\frac{ZZ^*}{|Z|^2}
\]
embeds $\mathbb P^{N-1}$ into the Euclidean space of Hermitian
matrices of trace one, with inner product
$\langle A,B\rangle=\operatorname{tr}(AB)$.
For $|Z|=1$ and $V\perp Z$,
\[
 dP_Z(V)=VZ^*+ZV^*,\qquad |dP_Z(V)|^2=2|V|^2.
\]
Thus $P$ is isometric for our Fubini--Study normalization.

At a point of $X$, choose a unitary basis $e_0,\ldots,e_{N-1}$
with $Z=e_0$ and $T^{1,0}X=\operatorname{span}_{\mathbb C}
\{e_1,\ldots,e_n\}$. Put $E_{ii}=e_ie_i^*$. The two
unit-speed projective geodesics associated with $e_i$ are
\[
 Z_i(s)=\cos(s/\sqrt2)e_0+\sin(s/\sqrt2)e_i,\qquad
 \widetilde Z_i(s)=\cos(s/\sqrt2)e_0+
            \sqrt{-1}\sin(s/\sqrt2)e_i.
\]
Their second derivatives under $P$ satisfy
\[
 (P\circ Z_i)''(0)=(P\circ\widetilde Z_i)''(0)
                    =E_{ii}-E_{00}.
\]
The complex submanifold $X$ has zero mean curvature in projective
space. Summing the second derivatives computed above over the
$2n$ real tangent directions therefore gives the Euclidean
mean-curvature vector
\[
 H=2\sum_{i=1}^nE_{ii}-2nE_{00},\qquad |H|^2=4n(n+1).
\]

Set
\[
 c=m\left(\frac{D|B^D|}{k|B^k|}\right)^{1/m},\qquad
 H_0=2\sqrt{n(n+1)},\qquad a_0=\frac{2(m-1)}{m-2}.
\]
Apply Lemma~\ref{lem:brendle} to the smooth positive function
$(w^2+\varepsilon^2)^{a_0/2}$ and let $\varepsilon\downarrow0$.
The Cauchy--Schwarz inequality then gives
\[
 c\|w\|_{2m/(m-2),dV}^{a_0}
 \le\|w\|_{2m/(m-2),dV}^{a_0-1}
       (a_0\|\nabla w\|_{2,dV}+H_0\|w\|_{2,dV}).
\]
Divide by $\|w\|_{2m/(m-2),dV}^{a_0-1}$ when $w\ne0$;
the resulting inequality also holds for $w=0$.

The real volume is $V=d(2\pi)^n/n!$, and $d\mu_0=dV/V$.
Since $a_0^2\le H_0^2$, a second application of the
Cauchy--Schwarz inequality yields
\[
 \|w\|_{2m/(m-2),\mu_0}^2
 \le S_*\bigl(\|\nabla w\|_{2,\mu_0}^2+\|w\|_{2,\mu_0}^2\bigr),
 \qquad S_*=2H_0^2c^{-2}V^{2/m}.
\]
Using
$D|B^D|/(k|B^k|)=\pi^n\Gamma(k/2)/\Gamma(k/2+n)$, we obtain
\[
 S_*=\frac{4(n+1)}n\left(\frac d{n!}
       \prod_{j=0}^{n-1}\left(\frac k2+j\right)\right)^{1/n}.
\]
For $n=2$, $S_*^2=2574<51^2$; for $n=3$,
$S_*^3=4096\cdot825/27<51^3$. Finally,
$|\nabla w|^2=2|\partial w|^2$ proves
\eqref{eq:appendix-projective-S}.
\end{proof}

\subsection{A Poincar\'e bound}
The proof below uses the logarithmic-gradient method of
Cheng and Yau \cite{ChengYau} and Li and Yau \cite{LiYau}.
We give the calculation for the constant needed here.
\begin{proposition}[A Ricci--diameter bound]\label{prop:ricci-diameter-gap}
Let $(M^m,g)$ be a connected closed Riemannian manifold with
$\operatorname{Ric}\ge-K_0g$, $K_0\ge0$, and diameter at most
$D_0>0$. Its first nonzero eigenvalue of $-\Delta$ satisfies
\begin{equation}
 \lambda_1\ge D_0^{-2}
 \exp\left[-\sqrt{(2m+2)(3+2K_0D_0^2)}\right].
\label{eq:appendix-geometric-gap}
\end{equation}
\end{proposition}
\begin{proof}
Normalize a first eigenfunction $v$ by $\max v=1$ and
$\min v=-a$, where $0<a\le1$, changing its sign if necessary.
For $\varepsilon>0$, set
\[
 u=v+a,\quad w=\log(u+\varepsilon),\quad F=|\nabla w|^2,
 \quad s=\lambda_1\frac{u-a}{u+\varepsilon},
 \quad L=\lambda_1(1+\varepsilon^{-1}).
\]
Then
\[
 \Delta w=-F-s,\qquad |s|\le L,\qquad
 \nabla s=\lambda_1\frac{a+\varepsilon}{u+\varepsilon}\nabla w,
 \qquad0\le\lambda_1\frac{a+\varepsilon}{u+\varepsilon}\le L.
\]
Apply the Bochner identity at a maximum point of $F$:
\[
 0\ge\frac2m(F+s)^2-2(L+K_0)F
 \ge\frac1mF^2-\frac2mL^2-2(L+K_0)F.
\]
Solving this inequality bounds the maximum, and hence gives on $M$:
\[
 F\le m(L+K_0)+\sqrt{m^2(L+K_0)^2+2L^2}
   \le(2m+2)(L+K_0).
\]
Integrate $|\nabla w|$ along a minimizing geodesic joining a
minimum and a maximum of $u$. Since their values are $0$ and
$1+a$, respectively,
\[
 \log(1+\varepsilon^{-1})
 \le D_0\sqrt{(2m+2)(L+K_0)}.
\]
If $\lambda_1>D_0^{-2}$, the proposition follows. Otherwise choose
\[
 \varepsilon=\frac{\lambda_1}{\lambda_1+K_0+D_0^{-2}}.
\]
Then $L+K_0\le3D_0^{-2}+2K_0$, and hence
\[
 \log\frac{D_0^{-2}}{\lambda_1}
 \le\log\left(2+\frac{K_0+D_0^{-2}}{\lambda_1}\right)
 \le\sqrt{(2m+2)(3+2K_0D_0^2)}.
\]
Exponentiation proves \eqref{eq:appendix-geometric-gap}.
\end{proof}

\subsection{Uniform constants on the family}
\begin{lemma}\label{lem:parametric-background}
For either parameter set $B_d$ in \eqref{eq:parameter-balls},
there is a rational $\lambda_0>0$ such that
\eqref{eq:background-P} and \eqref{eq:background-S} hold for every
$(X_t,\omega_{0,t})$, with $S_0=102$ and the same $\lambda_0$.
A finite search through rational inequalities for the defining
polynomials determines such a $\lambda_0$.
\end{lemma}
\begin{proof}
Proposition~\ref{prop:projective-sobolev} gives $S_0=102$.
To obtain a uniform Poincar\'e constant, we seek rational
numbers $K_0\ge0$ and $D_0>0$ with
\begin{equation}
 \operatorname{Ric}(g_{0,t})\ge-K_0g_{0,t},\qquad
 \operatorname{diam}(X_t,g_{0,t})<D_0\qquad(t\in B_d).
 \label{eq:family-background-geometry}
\end{equation}

In an affine patch, split the real coordinates into
$x\in\mathbb R^{2n}$ and $y\in\mathbb R^2$, where $y$ is one
complex coordinate, and put
$G=(\operatorname{Re}F_t,\operatorname{Im}F_t)$.
Choose closed boxes with rational endpoints
\[
 J_a\subset\operatorname{int}J_a^+,\qquad
 B\subset\operatorname{int}B^+,\qquad Q^-\subset\operatorname{int}Q.
\]
Require, for every $t\in J_a^+$,
\begin{align}
 &\forall x\in B^+\ \exists y\in Q^-:\ G(t,x,y)=0;\notag\\
 &x\in B^+,\ y,y'\in Q,\ G(t,x,y)=G(t,x,y')=0
                         \ \Longrightarrow\ y=y';\notag\\
 &x\in B^+,\ y\in Q,\ G(t,x,y)=0
              \ \Longrightarrow\ |\det D_yG(t,x,y)|\ge\eta,
 \label{eq:family-graph-conditions}
\end{align}
where $\eta>0$ is rational. The implicit-function theorem gives
a smooth graph $y=h(t,x)$ on a neighborhood of
$J_a^+\times B^+$. Its first derivatives satisfy
\[
 D_xh=-(D_yG)^{-1}D_xG,\qquad D_th=-(D_yG)^{-1}D_tG.
\]
Repeated differentiation expresses all derivatives needed for
the induced metric and its curvature as rational functions
of $(t,x,h(t,x))$, with nonzero denominators.

For each $J_a^+$, choose finitely many such graphs. Let $U_{ai,t}$
be the graph over the interior of its smaller box $B_{ai}$.
Require these graph neighborhoods to cover $X_t$. Choose a tree on
their indices such that
\[
 U_{ai,t}\cap U_{aj,t}\ne\varnothing
 \quad\text{for each edge }(i,j)\text{ and every }t\in J_a^+.
\]
Also require $B_d\subset\bigcup_a\operatorname{int}J_a$.
The cover and overlap requirements are polynomial conditions:
use homogeneous
representatives with $\sum|Z_i|^2=1$ and express each affine
coordinate by its defining polynomial equation.

On each enlarged graph, choose rational $b_{ai}>0$ and
$K_{ai}\ge0$ so that
\[
 g_{ai,t}\le b_{ai}I,\qquad
 \operatorname{Ric}(g_{0,t})(v,v)\ge-K_{ai}g_{0,t}(v,v)
       \quad(v\in TX_t).
\]
These inequalities are also polynomial conditions after
multiplication by the squares of the nonzero denominators.

A straight segment in $B_{ai}$ lifts to a path of length at most
$\sqrt{b_{ai}}\operatorname{diam}_{\mathbb R^{2n}}(B_{ai})$.
Choose a positive rational $d_{ai}$ with
\[
 d_{ai}^2>b_{ai}\operatorname{diam}_{\mathbb R^{2n}}(B_{ai})^2.
\]
Follow the chosen tree and join consecutive overlap points
within each graph. The segment bounds then give a path between
any two points of $X_t$ with length less than $\sum_i d_{ai}$. Thus
\[
 K_0=\max_{a,i}K_{ai},\qquad D_0=1+\max_a\sum_i d_{ai}
\]
satisfy \eqref{eq:family-background-geometry}.

To prove that the search stops, we first establish the existence
of the required finite data. On a fixed smooth connected fiber,
compactness gives a finite graph cover
with strict interior and Jacobian bounds. Its intersection graph
is connected, so it contains a spanning tree. Choose one point
in each overlap belonging to an edge of that tree.

The implicit-function theorem continues the graphs and the
chosen overlap points to nearby parameters. The cover also
persists: a sequence of uncovered points would, by properness,
have an uncovered limit on the original fiber.

Choose rational boxes and rational metric bounds within the
resulting strict margins. Compactness of $B_d$ then gives
finitely many parameter boxes.

Enumerate rational boxes, rational bounds, and finite trees in an
order with finitely many predecessors for each candidate.
For each candidate, Lemma~\ref{lem:real-elimination} decides
whether the required conditions hold. The preceding existence
argument supplies an admissible candidate, and the enumeration
reaches it after finitely many steps.

Finally, choose an integer $N>0$ with
\[
 N^2>(4n+2)(3+2K_0D_0^2),\qquad n=d-2,
\]
and put $\lambda_0=(2D_0^2\,4^N)^{-1}$.
Proposition~\ref{prop:ricci-diameter-gap} and
$\Delta_{\mathbb R}=2\Delta_{\omega_0}$ give
\[
 \frac{\lambda_1(-\Delta_{\mathbb R})}{2}
 >\frac1{2D_0^2e^N}>\lambda_0.
\]
The Rayleigh quotient therefore gives \eqref{eq:background-P}
with this rational constant.
\end{proof}

\section{Pointwise comparison}\label{sec:pointwise}
\subsection{An explicit bound for the potential}
\begin{proposition}[Potential bound]\label{prop:explicit-C0}
Let $(X,\omega)$ be a compact connected K\"ahler manifold of
complex dimension $n\in\{2,3\}$. Suppose
\eqref{eq:background-P} and \eqref{eq:background-S} hold with
constants $\lambda>0$ and $S\ge1$. If
\[
 \gamma=\omega+\ddbar v>0,\qquad
 \gamma^n=f\omega^n,\qquad
 \int_Xv\,d\mu_\omega=0,\qquad
 \delta=\|f-1\|_\infty,
\]
where $v$ and $f>0$ are smooth real-valued functions, then
\begin{align}
 \|v\|_2&\le n\delta/\lambda,\qquad
 \|\partial v\|_2\le n\delta/\sqrt\lambda,\label{eq:effective-potential-L2}\\
 \|v\|_\infty&\le C_\infty\delta,\qquad
 C_\infty=[S(n+1)]^{n/2}
       \left(\frac n{n-1}\right)^{n(n-1)/2}
       \max\{1,n/\lambda\}.\label{eq:explicit-C0}
\end{align}
All norms use $\omega$ and $d\mu_\omega=\omega^n/\int_X\omega^n$.
There is no restriction on the size of $\delta$.
\end{proposition}
\begin{proof}
Write $V=\int_X\omega^n$, $d\mu=d\mu_\omega$, and set
\[
 \Theta=
 \begin{cases}
  \gamma+\omega,&n=2,\\
  \gamma^2+\gamma\wedge\omega+\omega^2,&n=3.
 \end{cases}
\]
Then $\Theta$ is closed, $\Theta\ge\omega^{n-1}$, and
\begin{equation}
 \ddbar v\wedge\Theta=(f-1)\omega^n.
 \label{eq:potential-factorization}
\end{equation}
Multiply \eqref{eq:potential-factorization} by $v$ and integrate
by parts. The identity
$\sqrt{-1}\partial v\wedge\bar\partial v\wedge\omega^{n-1}
=n^{-1}|\partial v|_\omega^2\omega^n$ gives
\[
 \frac1n\|\partial v\|_2^2
 \le\frac1V\int_X\sqrt{-1}\partial v\wedge\bar\partial v\wedge\Theta
 =-\int_Xv(f-1)\,d\mu\le\delta\|v\|_2.
\]
Together with \eqref{eq:background-P}, this proves
\eqref{eq:effective-potential-L2}. It also gives $v=0$ if $\delta=0$.

For $\delta>0$, put $w=v/\delta$. Multiply
\eqref{eq:potential-factorization} by $|w|^{p-2}w/\delta$,
where $p\ge2$, and integrate by parts. This gives
\[
 \frac{p-1}{n}\int_X|w|^{p-2}|\partial w|_\omega^2\,d\mu
 \le\int_X|w|^{p-1}\,d\mu,
\]
and hence
\[
 \bigl\|\partial|w|^{p/2}\bigr\|_2^2
 \le\frac{np^2}{4(p-1)}\|w\|_{p-1}^{p-1}.
\]
Integration by parts is valid for the $C^1$ multiplier
$|w|^{p-2}w$; the Sobolev estimate applies to $|w|^{p/2}$.
Set $a_p=\max\{1,\|w\|_p\}$. Since $\int_Xd\mu=1$,
the gradient estimate and \eqref{eq:background-S} give
\begin{align*}
 \|w\|_{pn/(n-1)}^p
 &\le S\left(\bigl\|\partial|w|^{p/2}\bigr\|_2^2+\|w\|_p^p\right)\\
 &\le S\left(\frac{np^2}{4(p-1)}+1\right)a_p^p
 \le S(n+1)\frac p2\,a_p^p.
\end{align*}
The Moser iteration \cite{Moser}, with $p_k=2(n/(n-1))^k$, gives
\[
 a_{p_{k+1}}\le[S(n+1)p_k/2]^{1/p_k}a_{p_k},\qquad
 \sum_{k\ge0}\frac1{p_k}=\frac n2,\quad
 \sum_{k\ge0}\frac{k}{p_k}=\frac{n(n-1)}2.
\]
Letting $k\to\infty$ and using $a_2\le\max\{1,n/\lambda\}$
proves \eqref{eq:explicit-C0}.
\end{proof}

\subsection{Comparisons under a curvature bound}
Let $X\subset\mathbb P^{n+1}$ be a smooth quartic surface or
quintic threefold, and put $\omega_0=\omega_{\mathrm{FS}}|_X$.
Let $\beta$ be the Ricci-flat form in $[\omega_0]$. For a
K\"ahler form $\alpha\in[\omega_0]$, write
\[
 \beta=\alpha+\ddbar u,\qquad q=\frac{\beta^n}{\alpha^n},
 \qquad\int_Xu\,d\mu_\alpha=0,\qquad
 d\mu_\alpha=\frac{\alpha^n}{\int_X\alpha^n}.
\]

\begin{lemma}[Transfer of the background inequalities]
\label{lem:uniform-background-transfer}
Suppose \eqref{eq:background-P} and \eqref{eq:background-S} hold
for $\omega_0$, with constants $S_0\ge1$ and $\lambda_0>0$.
If $0<\ell\le1\le M$ and $\ell\omega_0\le\alpha\le M\omega_0$,
then the same inequalities hold for $\alpha$ with
\[
 S_\alpha=S_0\frac{M^{n-1}}{\ell^n},\qquad
 \lambda_\alpha=\lambda_0\frac{\ell^{n-1}}{M^n}.
\]
Define
\begin{equation}
 C_\alpha=[(n+1)S_\alpha]^{\lceil n/2\rceil}
 \left(\frac n{n-1}\right)^{n(n-1)/2}
 \max\{1,n/\lambda_\alpha\}.
 \label{eq:uniform-rational-C}
\end{equation}
Then $\operatorname{osc}u\le2C_\alpha\|q-1\|_\infty$.
The number $C_\alpha$ depends only on $n,S_0,\lambda_0,\ell,M$
and is rational if $S_0,\lambda_0,\ell,M$ are rational.
\end{lemma}
\begin{proof}
The forms $\alpha$ and $\omega_0$ have equal total volumes.
In an $\omega_0$-unitary frame diagonalizing $\alpha$, let
$a_i\in[\ell,M]$ be the eigenvalues. Then
\[
 d\mu_\alpha=D\,d\mu_0,\qquad D=\prod_i a_i,\qquad
 \ell^n\le D\le M^n,
\]
and
\[
 |\partial w|_\alpha^2D
 =\sum_i\left(\prod_{j\ne i}a_j\right)|w_i|^2
 \ge\ell^{n-1}|\partial w|_0^2.
\]
Consequently, every smooth real-valued function $w$ satisfies
\begin{align*}
 \|w\|_{2n/(n-1),\mu_\alpha}^2
 &\le M^{n-1}\|w\|_{2n/(n-1),\mu_0}^2\\
 &\le S_0M^{n-1}\left(\|\partial w\|_{2,\mu_0}^2+\|w\|_{2,\mu_0}^2\right)\\
 &\le\frac{S_0M^{n-1}}{\ell^n}
       \left(\ell\|\partial w\|_{2,\mu_\alpha}^2+\|w\|_{2,\mu_\alpha}^2\right).
\end{align*}
Since $\ell\le1$, this is the claimed Sobolev inequality.

To obtain the Poincar\'e inequality for $\alpha$, suppose
$\int_Xw\,d\mu_\alpha=0$ and set $\bar w_0=\int_Xw\,d\mu_0$.
The zero-mean condition, density comparison, and background
Poincar\'e inequality give
\begin{align*}
 \|w\|_{2,\mu_\alpha}^2
 &=\inf_{c\in\mathbb R}\int_X|w-c|^2\,d\mu_\alpha\\
 &\le M^n\int_X|w-\bar w_0|^2\,d\mu_0
 \le\frac{M^n}{\lambda_0}\|\partial w\|_{2,\mu_0}^2\\
 &\le\frac{M^n}{\lambda_0\ell^{n-1}}\|\partial w\|_{2,\mu_\alpha}^2.
\end{align*}
Finally, Proposition~\ref{prop:explicit-C0} gives
$\operatorname{osc}u\le2C_\infty\|q-1\|_\infty$.
As $(n+1)S_\alpha\ge1$, replacing $n/2$ by $\lceil n/2\rceil$
in $C_\infty$ gives the rational upper bound
\eqref{eq:uniform-rational-C}.
\end{proof}

\begin{proposition}[Pointwise comparison]
\label{prop:uniform-criterion}
Let $\alpha$ and $\beta$ be cohomologous K\"ahler forms on a
compact connected complex manifold $X$ of dimension $n\in\{2,3\}$,
with $\beta$ Ricci-flat. Let $\Omega$ be a nowhere vanishing
holomorphic $n$-form, and write
\[
 \beta=\alpha+\ddbar u,\qquad
 R=\frac{\alpha^n/n!}{(\sqrt{-1})^{n^2}\Omega\wedge\overline\Omega}.
\]
Suppose $0<L\le R\le U$ and
\[
 R^\alpha(\xi,\bar\xi,\zeta,\bar\zeta)
       \le K|\xi|_\alpha^2|\zeta|_\alpha^2,\qquad K\ge0.
\]
Let $b\ge\operatorname{osc}u$ be an oscillation bound and choose
a real number $a>K$ with $ab<1$. Set $\kappa=U/L$ and
\begin{equation}
 \tau=\frac{na}{(a-K)(1-ab)}.
 \label{eq:uniform-TE}
\end{equation}
The number $\tau$ bounds $\operatorname{tr}_\beta\alpha$.
If $0<r<1$ and $2\tau(\tau-n+\kappa-1)<r^2$, then
\begin{equation}
 \frac1{1+r}\alpha<\beta<\frac1{1-r}\alpha.
 \label{eq:uniform-comparison}
\end{equation}
In particular,
\begin{align}
 \|\beta-\alpha\|_{C^0(\alpha),\op}
   &<\frac r{1-r},\label{eq:uniform-relative-forward}\\
 \|\alpha-\beta\|_{C^0(\beta),\op}
   &<r.\label{eq:uniform-relative-reverse}
\end{align}
If $\alpha\le M\omega_0$, then
\begin{equation}
 \|\beta-\alpha\|_{C^0(\omega_0),\op}<\frac{Mr}{1-r}.
 \label{eq:uniform-fixed}
\end{equation}
If $\ell\omega_0\le\alpha\le M\omega_0$, with
$0<\ell\le1\le M$, and the background inequalities hold with
$S_0,\lambda_0$, one may choose
$b=2C_\alpha(\kappa-1)$, where $C_\alpha$ is given by
\eqref{eq:uniform-rational-C}.
\end{proposition}
\begin{proof}
Let $dV_\Omega=(\sqrt{-1})^{n^2}\Omega\wedge\overline\Omega$
and $d\nu=dV_\Omega/\int_XdV_\Omega$. Ricci-flatness and
equality of the K\"ahler classes give
\[
 \frac{\beta^n/n!}{dV_\Omega}=C,\qquad
 C=\int_XR\,d\nu\in[L,U],\qquad
 q:=\frac{\beta^n}{\alpha^n}=\frac C R.
\]
Therefore
\[
 \kappa^{-1}\le q\le\kappa,\qquad
 \|q-1\|_\infty\le\kappa-1.
\]
When the hypotheses of
Lemma~\ref{lem:uniform-background-transfer} hold, apply it together
with the density bound just proved. The stated choice of $b$ then bounds
$\operatorname{osc}u$.

Apply the Chern--Lu formula to
$\mathrm{id}:(X,\beta)\to(X,\alpha)$; see Yau
\cite[p.~199, (16)]{YauSchwarz}. In a $\beta$-unitary frame
$e_1,\ldots,e_n$, put $V_i=d\mathrm{id}(e_i)$ and
$s=\sum_i|V_i|_\alpha^2=\operatorname{tr}_\beta\alpha$.
With the curvature and Laplacian conventions of
Section~\ref{sec:preliminaries},
\[
 \Delta_\beta s
 =|\nabla d\mathrm{id}|^2
   -\sum_{i,k}R^\alpha(V_i,\overline{V_i},V_k,\overline{V_k})
 \ge|\nabla d\mathrm{id}|^2-Ks^2.
\]
The Cauchy--Schwarz inequality gives
$|\partial s|_\beta^2\le s|\nabla d\mathrm{id}|^2$. Hence
\[
 \Delta_\beta\log s\ge-Ks,\qquad
 \Delta_\beta u=n-s,\qquad
 \Delta_\beta(\log s-au)\ge(a-K)s-na.
\]
At a maximum point $x_0$ of $\log s-au$, the maximum
principle gives
$s(x_0)\le na/(a-K)$. The maximality then gives, for every $x\in X$,
\begin{equation}
 s(x)\le s(x_0)e^{a(u(x)-u(x_0))}
 \le\frac{na}{a-K}e^{ab}
 \le\frac{na}{(a-K)(1-ab)}=\tau.
 \label{eq:exponential-trace-bound}
\end{equation}
The last inequality uses $e^x\le(1-x)^{-1}$ for $0\le x<1$.

To bound the eigenvalues from the trace, put $h(t)=t-1-\log t$. Then
\[
 h(t)=
 \begin{cases}
  \displaystyle\int_t^1\frac{1-s}{s}\,ds
                    \ge\frac{(t-1)^2}{2},&0<t\le1,\\[6pt]
  \displaystyle\int_1^t\frac{s-1}{s}\,ds
                    \ge\frac{(t-1)^2}{2t},&t\ge1.
 \end{cases}
\]
Thus
\begin{equation}
 h(t)\ge\frac{(t-1)^2}{2\max\{1,t\}}\qquad(t>0).
 \label{eq:scalar-entropy-square}
\end{equation}

Let $\mu_i$ be the eigenvalues of $\beta^{-1}\alpha$.
Equations \eqref{eq:exponential-trace-bound} and
$\prod_i\mu_i=q^{-1}$ imply
\[
 \sum_i h(\mu_i)=s-n+\log q
 \le\tau-n+\log\kappa\le\tau-n+\kappa-1.
\]
Since $h(\mu_i)\ge0$ and $\mu_i\le s\le\tau$, with $\tau\ge n$,
\eqref{eq:scalar-entropy-square} gives
\[
 (\mu_i-1)^2\le2\tau h(\mu_i)
 \le2\tau(\tau-n+\kappa-1)<r^2.
\]
Hence $1-r<\mu_i<1+r$, which proves
\eqref{eq:uniform-comparison} and
\eqref{eq:uniform-relative-reverse}. Taking reciprocals gives
\eqref{eq:uniform-relative-forward}. Finally,
\[
 -\frac r{1+r}\alpha<\beta-\alpha<\frac r{1-r}\alpha,
\]
so $\alpha\le M\omega_0$ proves \eqref{eq:uniform-fixed}.
\end{proof}

\section{The construction algorithm}
\label{subsec:parametric-uniform}
We construct the polynomials in Theorem~\ref{thm:strong-intro}.
Throughout this section, $d=4$ or $5$, $n=d-2$, and $t\in B_d$,
with $B_d$ as in \eqref{eq:parameter-balls}.

\begin{lemma}[Polynomial approximation]\label{lem:density-potentials}
Let $K$ be a compact subset of a smooth embedded submanifold
$Y\subset\mathbb R^D$, and let $u$ be a smooth real-valued
function on a neighborhood of $K$ in $Y$. For each integer
$r\ge0$ and each $\varepsilon>0$, there is a polynomial $P$
in the ambient real coordinates, with rational coefficients,
such that $\|P|_Y-u\|_{C^r(K)}<\varepsilon$. Here the norm
is taken in a fixed finite collection of charts near $K$.
\end{lemma}
\begin{proof}
Extend $u$ to a compactly supported smooth real function $F$
on $\mathbb R^D$. This is possible by local tubular
neighborhoods, a partition of unity and a cutoff near $K$.
For $a>0$, put
\[
 F_a(x)=(\pi a^2)^{-D/2}
       \int_{\mathbb R^D}e^{-|x-y|^2/a^2}F(y)\,dy.
\]
Then $F_a\to F$ in $C^r$ on every compact set as $a\to0$.
The function $F_a$ extends to an entire function on $\mathbb C^D$
by replacing $|x-y|^2$ in this formula with
$\sum_{i=1}^D(z_i-y_i)^2$. Its Taylor polynomials therefore
converge, with their derivatives through order $r$, on a
closed box containing $K$. These polynomials have real coefficients.

Choose one Taylor polynomial sufficiently close to $F$ and
replace its finitely many coefficients by nearby rational numbers.
This changes its $C^r$ norm on the box by an arbitrarily small
amount. Restriction to the fixed charts on $Y$ proves the claim.
\end{proof}

\begin{lemma}[The Calabi--Yau potentials in the family]
\label{lem:family-CY-potentials}
There is a smooth real-valued function $\psi$ on a neighborhood
of the total family over $B_d$ such that
\[
 \beta_t=\omega_{0,t}+\sqrt{-1}\partial_{X_t}\bar\partial_{X_t}\psi_t,
 \qquad
 \int_{X_t}\psi_t\,d\mu_{0,t}=0.
\]
For each integer $k\ge0$, polynomials in the real and imaginary
parts of $(t,\Pi)$, with rational coefficients, approximate
$\psi$ in $C^k$ on this compact total family.
\end{lemma}
\begin{proof}
The Calabi--Yau theorem gives the normalized potential on each
fiber. By Lemma~\ref{lem:family-smoothness}, these fibers form
a smooth proper family over a neighborhood of $B_d$.
In a differentiable local trivialization, consider
\[
 (v,c)\longmapsto
 \left(
 \log\frac{(\omega_{0,t}+\sqrt{-1}\partial\bar\partial v)^n/n!}
                 {dV_{\Omega_t}}-c,
 \int_{X_t}v\,d\mu_{0,t}
 \right)
\]
as a map from $C^{k+2,a}\times\mathbb R$ to
$C^{k,a}\times\mathbb R$, where $0<a<1$. At the
Calabi--Yau solution its derivative is
\begin{equation}
 (w,b)\longmapsto
 \left(\Delta_{\beta_t}w-b,
               \int_{X_t}w\,d\mu_{0,t}\right).
 \label{eq:normalized-linearization}
\end{equation}
This map is invertible. Given $(h,\tau)$, its inverse is
determined by
\[
 b=-\int_{X_t}h\,d\nu_t,
 \qquad \Delta_{\beta_t}w=h+b,
 \qquad \int_{X_t}w\,d\mu_{0,t}=\tau.
\]
The first equality makes the right-hand side of the Laplace
equation have zero integral against $d\nu_t$; the last fixes
the additive constant. Schauder estimates give a bounded
inverse; see, e.g., Gilbarg and Trudinger \cite[Theorem~6.2 and Corollary~6.3]{GT}.

The implicit-function theorem gives smooth dependence on $t$
in every $C^{k+2,a}$ norm. By uniqueness, the normalized
solutions agree on overlapping parameter neighborhoods, so
$\psi$ is smooth on the total family. Applying
Lemma~\ref{lem:density-potentials} to the embedding
$(t,[Z])\mapsto(t,\Pi(Z))$ gives the required polynomial approximations.
\end{proof}

\begin{proposition}[Selection by rational inequalities]
\label{prop:finite-polynomial-selection}
For each integer $j\ge1$, one can find, after finitely many
rational polynomial tests, a polynomial $P_j(t,\Pi)$ with
rational coefficients such that $\alpha_{j,t}>0$,
\eqref{eq:intro-effective-uniform} holds, and
\begin{equation}
 \begin{gathered}
 \sup_{t\in B_d}\max_{\gamma_t\in\{\alpha_{j,t},\beta_t\}}
 \|\alpha_{j,t}-\beta_t\|_{C^0(\gamma_t),\op}<2^{-j},\\
 \sup_{t\in B_d}
 \frac{\sup_{X_t}(\alpha_{j,t}^n/\beta_t^n)}
      {\inf_{X_t}(\alpha_{j,t}^n/\beta_t^n)}<1+2^{-j},\\
 \sup_{t\in B_d}\sup_{X_t}
 |\nabla_{0,t}^{\,k}\log R_{P_j}|_{g_{0,t}}<2^{-j}
 \quad(1\le k\le j+2).
 \end{gathered}
 \label{eq:relative-approximation}
\end{equation}
\end{proposition}
\begin{proof}
Lemma~\ref{lem:parametric-background} supplies rational
constants $S_0=102$ and $\lambda_0>0$ for all the backgrounds
$(X_t,\omega_{0,t})$. Fix $j$ and write $\varepsilon=2^{-j}$.
We first give sufficient conditions on a polynomial $P$ and
then prove that a finite search finds one satisfying them.

Choose rational scalar constants
\begin{equation}
 0<\ell\le1\le M,\qquad K\ge0,\qquad a>K,
 \qquad 1<\kappa<1+\varepsilon.
 \label{eq:approximation-data}
\end{equation}
Define $C_\alpha$ by \eqref{eq:uniform-rational-C}, and set
\[
 b=2C_\alpha(\kappa-1),\qquad
 r=\frac{\varepsilon}{2M+\varepsilon}.
\]
Require the metric and curvature bounds
\eqref{eq:parametric-geometric-tests}, the density-ratio bound
\eqref{eq:parametric-density-test}, and the derivative bounds
\eqref{eq:log-density-derivative-tests} with $s=j+2$ and
$\eta=\varepsilon$. Also require
\begin{equation}
 ab<1,\qquad \tau=\frac{na}{(a-K)(1-ab)},\qquad
 2\tau(\tau-n+\kappa-1)<r^2.
 \label{eq:approximation-conditions}
\end{equation}
Here $b$ is the oscillation bound from
Lemma~\ref{lem:uniform-background-transfer}, $a$ is the chosen
parameter in Proposition~\ref{prop:uniform-criterion}, and
$\tau$ is its trace bound. All these numbers are rational.

These conditions imply the required pointwise bounds. Indeed,
\eqref{eq:normalized-density} gives
\[
 \kappa^{-1}\le\frac{\alpha_{P,t}^n}{\beta_t^n}
    =\frac{R_P}{C_t}\le\kappa,
 \qquad C_t=\int_{X_t}R_P\,d\nu_t.
\]
Proposition~\ref{prop:uniform-criterion} applies, and
\[
 \frac{Mr}{1-r}=\frac\varepsilon2,\qquad
 \frac r{1-r}=\frac\varepsilon{2M}\le\frac\varepsilon2,
 \qquad r<\frac\varepsilon2.
\]
Thus the metric difference is less than $\varepsilon$ in the
three operator norms based on $\omega_{0,t}$, $\alpha_{P,t}$
and $\beta_t$, uniformly for $t\in B_d$. The ratio of the
largest to the smallest value of $\alpha_{P,t}^n/\beta_t^n$
on $X_t$ is at most $\kappa<1+\varepsilon$.

For a smooth real function $v$ and a $g_{0,t}$-unit tangent
vector $V$,
\[
 (\sqrt{-1}\partial\bar\partial v)(V,J_tV)
 =\frac12\bigl(\nabla_{0,t}^2v(V,V)
                   +\nabla_{0,t}^2v(J_tV,J_tV)\bigr).
\]
Consequently, \eqref{eq:ricci-residue-density} gives
\[
 \|\operatorname{Ric}(\alpha_{P,t})\|_{C^0(\omega_{0,t}),\op}
 \le\|\nabla_{0,t}^2\log R_P\|_{C^0(g_{0,t})}<\varepsilon.
\]

It remains to prove that the conditions can be satisfied.
By Lemma~\ref{lem:family-CY-potentials}, there are polynomials
$Q_m\to\psi$ in $C^{j+4}$ on the compact total family. Hence
$\alpha_{Q_m,t}\to\beta_t$ uniformly in $C^{j+2}$. Choose
fixed rational $\ell,M,K$ for which the metric and curvature
bounds for $\beta_t$ are strict. These bounds continue to hold
for all sufficiently large $m$. Since $R_{\beta_t}=C_t>0$ is
constant on each fiber,
\begin{equation}
 \begin{gathered}
 \sup_{t\in B_d}
 \frac{\sup_{X_t}R_{Q_m}}{\inf_{X_t}R_{Q_m}}\longrightarrow1,\\
 \sup_{t\in B_d}\sup_{X_t}
 |\nabla_{0,t}^{\,k}\log R_{Q_m}|_{g_{0,t}}\longrightarrow0
 \quad(1\le k\le j+2).
 \end{gathered}
 \label{eq:existence-derivative-vanishing}
\end{equation}
Choose rational $\kappa_m>
\sup_t(\sup_{X_t}R_{Q_m}/\inf_{X_t}R_{Q_m})$ with
$\kappa_m\to1$. The constants $C_\alpha$ and $r$ remain fixed.
Choose one rational $a>K$ sufficiently large that
\[
 \tau_*:=\frac{na}{a-K},\qquad 2\tau_*(\tau_*-n)<r^2.
\]
As $m\to\infty$,
\[
 b_m=2C_\alpha(\kappa_m-1)\longrightarrow0,\qquad
 \tau_m=\frac{na}{(a-K)(1-ab_m)}\longrightarrow\tau_*.
\]
Consequently,
\[
 2\tau_m(\tau_m-n+\kappa_m-1)
 \longrightarrow2\tau_*(\tau_*-n)<r^2.
\]
Equations \eqref{eq:approximation-data} and
\eqref{eq:approximation-conditions} therefore hold for large $m$;
\eqref{eq:existence-derivative-vanishing} supplies all the
derivative bounds.

Enumerate the rational coefficients and exponents of $P$,
together with the rational numbers in \eqref{eq:approximation-data},
first by the length of their finite descriptions and then
lexicographically. Each candidate has finitely many predecessors.
Lemma~\ref{lem:QE} decides the metric, curvature, density and
derivative conditions, while \eqref{eq:approximation-conditions}
requires only rational comparisons. Select the first candidate
satisfying every condition, and denote its polynomial by $P_j$.
The existence argument above guarantees that the search stops
after finitely many tests.
\end{proof}

\begin{proposition}[Smooth convergence of the selected metrics]
\label{prop:selected-smooth-convergence}
For the polynomials $P_j$ selected in
Proposition~\ref{prop:finite-polynomial-selection}, and for
each integer $r\ge0$, there are constants $C_{d,r}<\infty$
and $J_{d,r}\ge1$ such that
\[
 \sup_{t\in B_d}\|\alpha_{j,t}-\beta_t\|_{C^r(g_{0,t})}
 \le C_{d,r}2^{-j}\qquad(j\ge J_{d,r}).
\]
\end{proposition}
\begin{proof}
Put
\begin{equation}
 F_{j,t}=\log\frac{\alpha_{j,t}^n}{\beta_t^n}
         =\log R_{P_j}-\log C_t.
 \label{eq:selected-log-density}
\end{equation}
The density-ratio bound and \eqref{eq:normalized-density} give
$\|F_{j,t}\|_{C^0}<\log(1+2^{-j})<2^{-j}$. The derivative
bounds in \eqref{eq:relative-approximation} give
\begin{equation}
 \sup_{t\in B_d}\|F_{j,t}\|_{C^k(g_{0,t})}\le2^{-j}
 \qquad(0\le k\le j+2).
 \label{eq:selected-density-smooth}
\end{equation}
Fix finitely many smooth family charts near the compact total
family over $B_d$. For each fixed $r$, conversion from covariant
to coordinate derivatives and the $C^{r+1}$ bound in
\eqref{eq:selected-density-smooth} imply
\begin{equation}
 \sup_{t\in B_d}\|F_{j,t}\|_{C^{r,1/2}}
 \le a_{d,r}2^{-j}\qquad(j\ge\max\{1,r-1\}).
 \label{eq:selected-density-holder-rate}
\end{equation}
Here $a_{d,r}$ depends only on $d,r$ and these charts.

Apply the implicit-function theorem at $(v,b)=(0,0)$ to
\[
 (v,b)\longmapsto
 \left(\log\frac{(\beta_t+\sqrt{-1}\partial\bar\partial v)^n}
                         {\beta_t^n}-b,
       \int_{X_t}v\,d\mu_{0,t}\right).
\]
The derivative is the isomorphism
\eqref{eq:normalized-linearization}. By compactness of $B_d$,
there are constants $\rho_{d,r}>0$ and $b_{d,r}<\infty$ such that,
if $\sup_t\|F_{j,t}\|_{C^{r,1/2}}<\rho_{d,r}$, there are
$v_{j,t}$ and real constants $c_{j,t}$ satisfying
\begin{align*}
 (\beta_t+\sqrt{-1}\partial\bar\partial v_{j,t})^n
     &=e^{F_{j,t}+c_{j,t}}\beta_t^n,\\
 \int_{X_t}v_{j,t}\,d\mu_{0,t}&=0,\\
 \sup_t\bigl(\|v_{j,t}\|_{C^{r+2,1/2}}+|c_{j,t}|\bigr)
     &\le b_{d,r}\sup_t\|F_{j,t}\|_{C^{r,1/2}}.
 \end{align*}
The resulting forms are positive. The constants
$\rho_{d,r}$ and $b_{d,r}$ depend only on $d,r$ and the fixed
charts.

Choose $J_{d,r}\ge\max\{1,r-1\}$ with
$a_{d,r}2^{-J_{d,r}}<\rho_{d,r}$. For $j\ge J_{d,r}$,
\begin{equation}
 \sup_t\|v_{j,t}\|_{C^{r+2,1/2}}
 \le b_{d,r}a_{d,r}2^{-j}.
 \label{eq:selected-potential-rate}
\end{equation}

Integration of the Monge--Amp\`ere equation gives
\[
 \int_{X_t}\beta_t^n
 =e^{c_{j,t}}\int_{X_t}e^{F_{j,t}}\beta_t^n
 =e^{c_{j,t}}\int_{X_t}\alpha_{j,t}^n
 =e^{c_{j,t}}\int_{X_t}\beta_t^n.
\]
Thus $c_{j,t}=0$. Uniqueness in the Calabi--Yau theorem yields
$\alpha_{j,t}=\beta_t+\sqrt{-1}\partial\bar\partial v_{j,t}$.
Applying $\sqrt{-1}\partial\bar\partial$ to
\eqref{eq:selected-potential-rate} proves the proposition.
Fixing the family charts once for each $d$ makes the final
constants depend only on $d,r$.
\end{proof}

\begin{proof}[Proof of Theorem~\ref{thm:strong-intro}]
For each $j$, choose $P_j$ by
Proposition~\ref{prop:finite-polynomial-selection}. By that
proposition, the forms $\alpha_{j,t}$ are K\"ahler and satisfy
\eqref{eq:intro-effective-uniform}.
Proposition~\ref{prop:selected-smooth-convergence} gives
\eqref{eq:intro-smooth-approximation} for this sequence.
The choice of $P_j$ is independent of $r$.
\end{proof}

\begin{corollary}[Algebraic quartics and quintics]
\label{cor:algebraic-hypersurface-algorithm}
Let $F$ have algebraic complex coefficients and define a smooth
degree-$d$ hypersurface $X=\{F=0\}\subset\mathbb P^{d-1}$,
where $d=4$ or $5$. The conclusions of
Theorem~\ref{thm:strong-intro} hold on this single hypersurface:
the algorithm takes the exact coefficients of $F$ and $j$ and
returns $P_j(\Pi)$ with rational coefficients. The constants
in the $C^r$ estimate depend on $F$ and $r$.
\end{corollary}
\begin{proof}
Represent each real and imaginary coefficient by its minimal
polynomial and a rational isolating interval. By adjoining
these isolated numbers as variables, polynomial tests over
the algebraic coefficients reduce to rational tests.

Proposition~\ref{prop:projective-sobolev} supplies $S_0=102$.
Applying the graph construction of Lemma~\ref{lem:parametric-background}
to this single smooth connected hypersurface gives
a rational $\lambda_0>0$. Its normalized Calabi--Yau potential
has rational polynomial approximations by
Lemma~\ref{lem:density-potentials}. The proof of
Proposition~\ref{prop:finite-polynomial-selection} therefore
applies with $t$ omitted, and
Proposition~\ref{prop:selected-smooth-convergence} gives the
stated smooth convergence.
\end{proof}

\begin{corollary}[Invariant potentials on the Fermat hypersurfaces]
\label{thm:fermat-invariant-uniform}
On $X_d$, the construction in
Proposition~\ref{prop:finite-polynomial-selection} can be carried
out with $G_d$-invariant polynomials $P_j(\Pi)$ with rational
coefficients. Their K\"ahler forms satisfy
\eqref{eq:relative-approximation} and the pointwise and smooth
estimates of Theorem~\ref{thm:strong-intro}, with $t=0$.
\end{corollary}
\begin{proof}
Uniqueness in the Calabi--Yau theorem makes the normalized
potential $\psi_d$ invariant under $G_d$. We show that averaging
over $G_d$ preserves rational coefficients.

Write the independent real coordinates of a Hermitian matrix as
$W_{ii}=x_{ii}$ and $W_{ij}=x_{ij}+\sqrt{-1}y_{ij}$ for $i<j$,
with $W_{ji}=\overline{W}_{ij}$. A rational polynomial in
$(x,y)$ becomes a polynomial in $W$ over $\mathbb Q(\sqrt{-1})$.
For a monomial $W^a=\prod_{i,j}W_{ij}^{a_{ij}}$, put
$w_i(a)=\sum_j a_{ij}-\sum_j a_{ji}$. Averaging the diagonal
root-of-unity action gives
\begin{equation}
 \frac1{d^d}\sum_{\zeta_0^d=\cdots=\zeta_{d-1}^d=1}
       \prod_i\zeta_i^{w_i(a)}W^a
 =\begin{cases}
 W^a,& d\mid w_i(a)\text{ for all }i,\\
 0,&\text{otherwise}.
 \end{cases}
 \label{eq:fermat-phase-average}
\end{equation}
Next average over permutations of the indices. Both averaging
operations preserve coefficients in $\mathbb Q(\sqrt{-1})$
and commute with complex conjugation on Hermitian matrices.
After returning to $(x,y)$, the coefficients are real and
belong to $\mathbb Q(\sqrt{-1})$, so they are rational.

Let $\mathcal A_dp=|G_d|^{-1}\sum_{g\in G_d}p\circ g$.
Each $g$ is an isometry of $g_0$, so, for every $k$,
\[
 \|\mathcal A_dp-\psi_d\|_{C^k(g_0)}
 \le\|p-\psi_d\|_{C^k(g_0)}.
\]
It follows that invariant rational polynomials approximate $\psi_d$
in every fixed $C^k$ norm. In the search of
Proposition~\ref{prop:finite-polynomial-selection} at $t=0$,
replace each candidate $p$ by $\mathcal A_dp$. Formula
\eqref{eq:fermat-phase-average} makes each averaging step a
rational computation. The approximation just proved guarantees
that the search stops. Propositions~\ref{prop:finite-polynomial-selection}
and \ref{prop:selected-smooth-convergence} give the conclusions.
\end{proof}

\enlargethispage{3\baselineskip}
\begin{remark}
The invariant polynomials use both diagonal and off-diagonal
entries of $\Pi$. No density assertion for polynomials in the
diagonal entries alone is needed.
\end{remark}

\clearpage
\appendix
\section{Coefficients of the explicit potentials}
\label{app:coefficients}
Tables~\ref{tab:recovered-Q12} and~\ref{tab:recovered-T12}
give the integers in \eqref{eq:recovered-Q12-T12}.
In each table the right half continues the left half.
Every unlisted coefficient is zero.

\begin{table}[!ht]
\centering
\caption{The coefficients of $P_{Q_{12}}$ in \eqref{eq:recovered-Q12-T12}.}
\label{tab:recovered-Q12}
\begin{tabular}{@{}rr@{\qquad}rr@{}}
\toprule
$b$ & $a_b$ & $b$ & $a_b$\\
\midrule
$(1,0,0)$ & $34386793216192$ & $(0,2,1)$ & $8807230326762496$\\
$(0,1,0)$ & $-40636570414656$ & $(1,0,2)$ & $4365198837150016$\\
$(0,0,1)$ & $32411009297195$ & $(2,2,0)$ & $-1805302581301248$\\
$(2,0,0)$ & $-256998055566080$ & $(3,0,1)$ & $-2155909931473280$\\
$(1,1,0)$ & $570792192738368$ & $(5,0,0)$ & $850997550637888$\\
$(0,2,0)$ & $-398461877041728$ & $(0,1,2)$ & $-11369198410431488$\\
$(1,0,1)$ & $-467739236371136$ & $(1,3,0)$ & $-1552218252518464$\\
$(3,0,0)$ & $887297975176256$ & $(2,1,1)$ & $5177299591634112$\\
$(0,1,1)$ & $840961842615936$ & $(4,1,0)$ & $-2016951335537152$\\
$(2,1,0)$ & $-2539829380646912$ & $(0,0,3)$ & $5005014538551488$\\
$(0,0,2)$ & $-442963018459904$ & $(0,4,0)$ & $566237808445888$\\
$(1,2,0)$ & $3173708343785984$ & $(1,2,1)$ & $9106908696448$\\
$(2,0,1)$ & $2105950681040320$ & $(2,0,2)$ & $-2828866047685760$\\
$(4,0,0)$ & $-1355783961929792$ & $(3,2,0)$ & $954695778505280$\\
$(0,3,0)$ & $-2085640122495872$ & $(4,0,1)$ & $346971377742016$\\
$(1,1,1)$ & $-7553623649681728$ & $(6,0,0)$ & $5806629677120$\\
$(3,1,0)$ & $3372690510374912$ &  & \\
\bottomrule
\end{tabular}
\end{table}

\clearpage
\begin{table}[!ht]
\centering
\caption{The coefficients of $P_{T_{12}}$ in \eqref{eq:recovered-Q12-T12}.}
\label{tab:recovered-T12}
\begin{tabular}{@{}rr@{\qquad}rr@{}}
\toprule
$b$ & $a_b$ & $b$ & $a_b$\\
\midrule
$(1,0,0,0)$ & $-132720746868375$ & $(0,2,1,0)$ & $183737466309009408$\\
$(0,1,0,0)$ & $-221250121729640$ & $(1,0,2,0)$ & $358726975800720384$\\
$(0,0,1,0)$ & $655767598543524$ & $(1,1,0,1)$ & $278375131866067200$\\
$(2,0,0,0)$ & $846481914019584$ & $(2,2,0,0)$ & $4062387225131520$\\
$(0,0,0,1)$ & $71394992506368$ & $(3,0,1,0)$ & $-33891287498503680$\\
$(1,1,0,0)$ & $7356707367707136$ & $(5,0,0,0)$ & $-15705435300613632$\\
$(0,2,0,0)$ & $-22092625927020288$ & $(0,1,2,0)$ & $-125645789091225600$\\
$(1,0,1,0)$ & $-19190595086008320$ & $(0,2,0,1)$ & $-172008216118121472$\\
$(3,0,0,0)$ & $-4872895557752832$ & $(1,0,1,1)$ & $-466321142003490048$\\
$(0,1,1,0)$ & $91763249817420288$ & $(1,3,0,0)$ & $-19537529563736832$\\
$(1,0,0,1)$ & $6552668616301824$ & $(2,1,1,0)$ & $138181747994243328$\\
$(2,1,0,0)$ & $-22968251723748096$ & $(3,0,0,1)$ & $-75681879086690304$\\
$(0,0,2,0)$ & $-88006380151728384$ & $(4,1,0,0)$ & $61536419966050560$\\
$(0,1,0,1)$ & $-66963360838857984$ & $(0,0,3,0)$ & $-22177699730856192$\\
$(1,2,0,0)$ & $109408791488620032$ & $(0,1,1,1)$ & $203880961703369472$\\
$(2,0,1,0)$ & $74967847943035392$ & $(0,4,0,0)$ & $22084322161181952$\\
$(4,0,0,0)$ & $17307283480421376$ & $(1,0,0,2)$ & $100021556369005056$\\
$(0,0,1,1)$ & $124354183036561920$ & $(1,2,1,0)$ & $-50129936162724864$\\
$(0,3,0,0)$ & $-41963016014158080$ & $(2,0,2,0)$ & $-75213597684821760$\\
$(1,1,1,0)$ & $-435937208103599616$ & $(2,1,0,1)$ & $104002936940709120$\\
$(2,0,0,1)$ & $-10736969237708544$ & $(3,2,0,0)$ & $-69260826054403584$\\
$(3,1,0,0)$ & $-36699934311591168$ & $(4,0,1,0)$ & $29344281744084480$\\
$(0,0,0,2)$ & $-32508196536883968$ & $(6,0,0,0)$ & $-2801048013315072$\\
\bottomrule
\end{tabular}
\end{table}
\clearpage

\section*{Funding}
Kefeng Liu is partially supported by the National Key Research and
Development Program of China (grant~2024YFA1013201), the Natural Science
Foundation of Chongqing (grant~CSTB2024NSCQ-LZX0040), and the Special
Project of Chongqing Municipal Science and Technology Bureau
(grant~2025CCZ015).

Liangdi Zhang is supported by the National Natural Science Foundation
of China (grant~12601080), the Scientific Research Foundation of
Chongqing University of Technology (grant~2026ZDZ012), the Youth
Project of the Science and Technology Research Program of Chongqing
Education Commission of China (grant~KJQN202601108), and the Chongqing Natural Science Foundation Innovation and Development Joint Fund (grant~CSTB2025NSCQ-LZX0059).

\section*{Declaration of AI use}
ChatGPT assisted with preliminary mathematical arguments,
numerical experiments, and language editing. The authors take full
responsibility for the content of this article.

\section*{Data availability}
The polynomial coefficients and the procedures for the finite
computations are included in this article. No external datasets
are used.

\raggedbottom

\end{document}